\documentclass[twocolumn]{autart}    

\usepackage{graphicx}          
\usepackage{mathptmx} 
\usepackage{times} 
\usepackage{amsmath} 
\usepackage{amssymb}  
\usepackage{color}
\usepackage{epstopdf}
\usepackage{subfigure}
\usepackage{tikz}
\usetikzlibrary{matrix,arrows,shapes.geometric}
\usepackage[all,cmtip]{xy}

\newtheorem{theorem}{Theorem}

\newtheorem{remark}{Remark}
\newtheorem{example}{Example}

\def\a {\alpha}
\def\b {\beta}

\def\s {\sigma}
\def\w {\omega}

\def\g {\gamma}
\def\l {\lambda}
\def\E {\mathcal{E}}

\begin{document}

\begin{frontmatter}

\title{Free-Endpoint Optimal Control of Inhomogeneous Bilinear Ensemble Systems \thanksref{footnoteinfo}} 

\thanks[footnoteinfo]{
This work was supported in part by the National Science Foundation under the awards CMMI-1462796 and ECCS-1509342. The material in this paper was not presented at any conference meetings.} 

\author[SJ]{Shuo Wang}\ead{swang35@wustl.edu},    
\author[SJ]{Jr-Shin Li}\ead{jsli@wustl.edu}            

\address[SJ]{Electrical and Systems Engineering, Washington University, St. Louis, Missouri, 63130, USA} 

%


\begin{abstract}                          
	Optimal control of bilinear systems has been a well-studied subject in the areas of mathematical and computational optimal control. However, effective methods for solving emerging optimal control problems involving an ensemble of deterministic or stochastic bilinear systems are underdeveloped. These burgeoning problems arise in diverse applications from quantum control and molecular imaging to neuroscience.	In this work, we develop an iterative method to find optimal controls for an inhomogeneous bilinear ensemble system with free-endpoint conditions. The central idea is to represent the bilinear ensemble system at each iteration as a time-varying linear ensemble system, and then solve it in an iterative manner. We analyze convergence of the iterative procedure and discuss optimality of the convergent solutions. The method is directly applicable to solve the same class of optimal control problems involving a stochastic bilinear ensemble system driven by independent additive noise processes. We demonstrate the robustness and applicability of the developed iterative method through practical control designs in neuroscience and quantum control.
\end{abstract}

\end{frontmatter}

\section{Introduction}
Controlling a population system consisting of a large number of structurally identical dynamic units is an essential step that enables many cutting-edge applications in science and engineering. For example, in quantum science and technology, synchronization engineering, and circadian biology, a central control task is to design exogenous forcing that guides individual subsystems in the population or ensemble to behave in a desired or an optimal manner  \cite{Cory97,Brent06,Li_PRA06,Kiss:2007:1886-1889,Ching2013a,Schaettler2002}. Such optimal ``broadcast'' control designs are of theoretical and computational challenge because, in practice, only a single or sparsely distributed control signals are available to engineer individual or collective behavior of many or a continuum of dynamical systems \cite{Phelps14,Li_TAC12_QCP,Li_TAC13}.

There exist numerous numerical methods for solving optimal control problems of nonlinear systems \cite{Rao09}, many of which rely heavily on applying effective discretization schemes to discretize the system dynamics and then implementing numerical optimizations to solve the resulting nonlinear programs (NLPs) \cite{Gong06}. Canonical methods include direct and indirect shooting methods \cite{StoerBulirsch1980,Bulirsch1992} and spectral collocation methods such as the pseudospectral method \cite{Gong08}. Implementing these commonly-used computational methods to solve optimal control problems involving an ensemble system may encounter low efficiency, slow convergence, and instability issues. It is because each subsystem in the ensemble has an identical structure so that the resulting discretized large-scale NLPs are equipped with a distinctive sparse structure, and, furthermore, each subsystem shares a common control input so that these NLPs involve highly localized and restrictive constraints \cite{Li_PNAS11}.

In this paper, we study the optimal control of a bilinear ensemble system with inhomogeneous natural and translational dynamics, which models a wide range of practical optimal control design problems across disciplines, for example, optimal pulse design in quantum control \cite{Daoyi2014,Li_JCP11} and molecular imaging \cite{Woods06}, motion planning of robots in the presence of uncertainty \cite{Becker12}, and optimal stimulation of spiking neurons \cite{Ching2013a}. In our previous study, we investigated ensemble control designs in these compelling applications and, in particular, focused on devising minimum-energy controls with fixed-endpoint constraints \cite{Li_CDC15}. 
Here, we consider 
free-endpoint quadratic optimal bilinear ensemble control problems and develop an iterative method to solve this class of problems without the use of 
numerical optimizations. The procedure is based on constructing and solving a corresponding optimal control problem involving a linear ensemble system at each iteration, which is numerically tractable as shown in our previous work \cite{Li_TAC11,Li_ACC12_SVD}. Moreover, the established iterative method is directly applicable to find optimal controls for stochastic bilinear systems driven by additive noise, such as Poisson counters and Brownian motion. We note that iterative methods have been developed for solving free-endpoint optimal control problems \cite{Hofer1988} or optimal tracking \cite{Cimen2004} of a single deterministic bilinear system. 
These pervious studies lay the foundation of our new developments towards solving 
optimal control problems involving a bilinear ensemble system governed by inhomogeneous drift and translational dynamics.

This paper is organized as follows. In the next section, we formulate the optimal control problem involving a single inhomogeneous bilinear system. We present the iterative method to solve this optimal control problem and show the convergence of the iterative algorithm by using the fixed-point theorem. In Section \ref{sec:ensemble}, we extend the iterative method to deal with optimal control problems for bilinear ensemble systems. In Section \ref{sec:example}, we illustrate the robustness and applicability of the iterative method through the examples of controlling spiking neurons in the presence of jump processes and pulse design in protein nuclear magnetic resonance (NMR) spectroscopy.

\section{Optimal Control of Inhomogeneous Bilinear Systems}
\label{sec:inhomogeneous}
In this paper, we study optimal control problems involving an ensemble of inhomogeneous bilinear systems with state-invariant drift and translational dynamics of the form
\begin{align} 
	\frac{d}{dt}{X(t,\b)} =& A(\b)X(t,\b)+B(\b)u(t) \nonumber \\ 
	\label{eq:bilinear_ensemble}
	&+\Big[\sum_{i=1}^m u_i(t) B_i(\b)\Big] X(t,\b)+ g(\b),
\end{align}
where $X(t,\b)\in\mathbb{R}^n$ denotes the state, $\b\in K\subset\mathbb{R}^d$ is the system parameter varying on the compact set $K$ in the $d$-dimensional Euclidean space; $u=(u_1,\ldots,u_m)^T$ is the control function with $u_i:[0,t_f]\to\mathbb{R}$ being piecewise continuous;
$A\in C(K;\mathbb{R}^{n\times n})$, $B\in C(K;\mathbb{R}^{n\times m})$, and $B_i\in C(K;\mathbb{R}^{n\times n})$, $i=1,\ldots,m$, are real matrices whose elements are continuous functions over $K$, and $g\in C(K; \mathbb{R}^n)$. Specifically, we consider the free-endpoint optimal control minimizing the cost functional involving the trade-off between the terminal cost and the energy of the control input. 

In the following, we develop an iterative procedure for solving this challenging optimal ensemble control problem. To fix the idea, we first illustrate the framework through a free-endpoint, finite-time, quadratic optimal control problem involving a single deterministic time-invariant inhomogeneous bilinear system. Namely, we consider the problem
\begin{align}
	\min &\quad J = \frac{1}{2}\int_0^{t_f} u^T(t)Ru(t) \, dt + \Vert x(t_f) - x_d \Vert _2^2, \nonumber\\
	\label{eq:oc1}
	{\rm s.t.} &\quad \dot{x}=Ax+Bu+\Big[\sum_{i=1}^m u_i B_i\Big]x + g , \tag{P1} 
\end{align}
where $x(t)\in\mathbb{R}^n$ is the state and $u(t)=(u_1,\ldots,u_m)^T\in\mathbb{R}^m$ is the control with each $u_i$ piecewise continuous; $A\in\mathbb{R}^{n\times n}$, $B_i\in\mathbb{R}^{n\times n}$, and $B\in\mathbb{R}^{n\times m}$ are constant matrices, and $g\in \mathbb{R}^{n}$ is a constant vector. In the cost functional, $R\in\mathbb{R}^{m\times m}\succ 0$ is a positive definite weight matrix for the control energy and $\Vert x(t_f) - x_d \Vert _2^2 = (x(t_f) - x_d)^T (x(t_f) - x_d)$ represents the terminal cost with respect to the desired state $x_d\in\mathbb{R}^n$. In addition, we can represent the time-invariant bilinear system in \eqref{eq:oc1} as $\dot{x}=Ax+Bu+\big[\sum_{j=1}^n x_j(t) N_j\big]u + g$, in which we write the bilinear term $\big(\sum_{i=1}^m u_iB_i\big)x=\big(\sum_{j=1}^n x_jN_j\big) u$ with $x_j$ being the $j^{th}$ element of $x$ and $N_j\in\mathbb{R}^{n\times m}$ for $j=i,\dots,n$. 

We now solve the optimal control problem \eqref{eq:oc1} by Pontryagin's maximum principle. The Hamiltonian of this problem is
$$H(x,u,p)=\frac{1}{2}u^TRu+p^T \big\{Ax+\big[B+(\sum_{j=1}^n x_jN_j)\big] u + g \big\},$$ 
where $p(t)\in\mathbb{R}^n$ is the co-state vector. The optimal control is then obtained by the necessary condition, $\frac{\partial H}{\partial u}=0$ (since the control $u$ is unconstrained), given by 
\begin{equation}
	\label{eq:u}
	u^* = -R^{-1}\Lambda^T p,
\end{equation}
where $\Lambda=B+\sum_{j=1}^n x_jN_j$, and the optimal trajectories of the state $x$ and the co-state $p$ satisfy, for $t\in[0,t_f]$,
\begin{align}
	\label{eq:state}
	\dot{x}_i & = \big[Ax\big]_i-\Big[\Lambda R^{-1}\Lambda^T p\Big]_i + g, \\
	\label{eq:costate}
	\dot{p}_i & = -\big[A^Tp\big]_i+p^T\Big[N_iR^{-1}\Lambda^T +\Lambda R^{-1} N_i^T\Big]p,
\end{align}
with the boundary conditions $x(0)=x_0$ and $p(t_f)=2(x(t_f)-x_d)$ from the transversality condition, where $x_i$, $p_i$ and $[\,\cdot\,]_i$, $i=1,\dots,n$, are the $i^{th}$ component of the respective vectors. By the following change of variables,
\begin{align}
	\label{eq:A}
	& \tilde{A}_{ij} = A_{ij}-\Big[(N_j R^{-1} \Lambda^T  +\Lambda R^{-1} N_j^T\big)p\Big]_i, \\
	\label{eq:B}
	& \tilde{B}R^{-1}\tilde{B}^T = BR^{-1}B^T-\big(\sum_{j=1}^n x_jN_j\big)R^{-1}\big(\sum_{j=1}^n x_jN_j\big)^T,
\end{align}
we can rewrite \eqref{eq:state} and \eqref{eq:costate} into the form
\begin{align}
	\label{eq:state1}
	\dot{x} &= \tilde{A}x-\tilde{B}R^{-1}\tilde{B}^Tp + g, \qquad x(0)=x_0, \\
	\label{eq:costate1}
	\dot{p} &= -\tilde{A}^Tp, \qquad \qquad \qquad \quad p(t_f)=2(x(t_f)-x_d),
\end{align}
which are in the similar form, with an additional inhomogeneous term $g$, of the canonical equations that characterize the optimal trajectories of the LQR problem \cite{Moore2007}.

\subsection{Iteration Procedures}
\label{sec:iterative}
If the matrices $\tilde{A}$ and $\tilde{B}R^{-1}\tilde{B}^T$ in \eqref{eq:A} and \eqref{eq:B}, respectively, were known, then the two-point boundary value problem (TPBVP) described in \eqref{eq:state1} and \eqref{eq:costate1} can be solved numerically, e.g., by shooting methods \cite{Bulirsch1992}. However, they are state dependent, i.e., $\tilde{A}(x,p)$ and $\tilde{B}R^{-1}\tilde{B}^T(x)$, so that the TPBVP is intractable. To overcome this, we propose to solve it in an iterative manner by considering the iteration equations
\begin{align}
	\label{eq:x_k}
	& \dot{x}^{(k)} = \tilde{A}^{(k-1)}x^{(k)} - \tilde{B}^{(k-1)}R^{-1}(\tilde{B}^{(k-1)})^Tp^{(k)} + g, \\ 
	\label{eq:p_k}
	& \dot{p}^{(k)}= -(\tilde{A}^{(k-1)})^Tp^{(k)},
\end{align}
with the boundary conditions $x^{(k)}(0)=x_0$ and $p^{(k)}(t_f)=2(x^{(k)}(t_f)-x_d)$ for all iterations $k=0,1,2,\ldots$, where the matrices $\tilde{A}^{(k-1)}(x^{(k-1)},p^{(k-1)})$ and $\tilde{O}^{(k-1)}(x^{(k-1)})\doteq\tilde{B}^{(k-1)}R^{-1}(\tilde{B}^{(k-1)})^T(x^{(k-1)})$ are defined according to \eqref{eq:A} and \eqref{eq:B}, given by
\begin{align}
	\label{eq:Ak}
	& \tilde{A}_{ij}^{(k)} = A_{ij}-\Big[(N_j R^{-1} (\Lambda^{(k)})^T+\Lambda^{(k)}R^{-1} N_j)p^{(k)}\Big]_i, \\
	\label{eq:Bk}
	&\tilde{O}^{(k)} = BR^{-1}B^T-\big(\sum_{j=1}^n x_j^{(k)}N_j\big)R^{-1}\big(\sum_{j=1}^n x_j^{(k)}N_j\big)^T,
\end{align}
with $\Lambda^{(k)} = B+\sum_{j=1}^n x_j^{(k)} N_j$. In order to solve for \eqref{eq:x_k} and \eqref{eq:p_k}, we let
\begin{equation}
	\label{eq:p^k}
	p^{(k)}(t) = K^{(k)}(t)x^{(k)}(t)+s^{(k)}(t),
\end{equation}
where $K^{(k)}(t)\in\mathbb{R}^{n\times n}$ and $s^{(k)}(t)\in\mathbb{R}^n$. Substituting this into  \eqref{eq:p_k} and using \eqref{eq:x_k} yields the Riccati equation
\begin{align}
	\label{eq:K}
	\dot{K}^{(k)}&=-K^{(k)}\tilde{A}^{(k-1)}-(\tilde{A}^{(k-1)})^TK^{(k)} +K^{(k)}\tilde{O}^{(k)}K^{(k)},
\end{align}
and  
\begin{align}
	\label{eq:Sk}
	\dot{s}^{(k)}=& -\Big[(\tilde{A}^{(k-1)})^T-K^{(k)}\tilde{O}^{(k-1)}\Big]s^{(k)} - K^{(k)}g,
\end{align}
with the terminal conditions $K^{(k)}(t_f)=2I_n$, where $I_n$ denotes the $n\times n$ identity matrix, and $s^{(k)}(t_f)=-2x_d$. Solving $K^{(k)}$ and $s^{(k)}$ with \eqref{eq:K} and \eqref{eq:Sk} and using \eqref{eq:p^k}, $x^{(k)}$ and $p^{(k)}$ in \eqref{eq:x_k} and \eqref{eq:p_k} can be obtained.

We would like to emphasize that solving these iterative equations is equivalent to solving the following inhomogeneous LQR problem at each iteration $k=0,1,2,\ldots$,
\begin{align}
	\min &\quad J=\frac{1}{2}\int_0^{t_f} (u^{(k)})^T(t)Ru^{(k)}(t) \, dt + \Vert x^{(k)}(t_f) - x_d \Vert _2^2, \nonumber\\
	\label{eq:oc2}
	{\rm s.t.} &\quad \dot{x}^{(k)} =\tilde{A}^{(k-1)} x^{(k)}+\tilde{B}^{(k-1)} u^{(k)} + g \tag{P2}\\  
	&\quad x^{(k)}(0)=x_0.  \nonumber
\end{align}
If the iterative procedure is convergent, we expect that the optimal solution of \eqref{eq:oc2} is approaching that of \eqref{eq:oc1}. The convergence behavior and criteria of this iterative method are analyzed in the following section.

\subsection{Convergence of the Iterative Method}
\label{sec:convergence}
In this section, we analyze the convergence of the iterative algorithm and show the dependence of the convergence on the choice of the weight matrix $R$ and on the translational dynamics $g$. 
To facilitate the proof, we introduce the following mathematical tools. Considering the Banach spaces, $\mathcal{X}\doteq C([0,t_f];\, \mathbb{R}^n)$, $\mathcal{Y}\doteq C([0,t_f];\, \mathbb{R}^{n\times n})$, and $\mathcal{Z}\doteq C([0,t_f];\, \mathbb{R}^n)$ with the norms
\begin{align*}
	\| x \|_{\alpha} &= \sup_{t\in[0,t_f]} \big[\|x(t)\|\exp(-\a t)\big], & \text{for}\ \ x\in\mathcal{X}, \\
	\| y \|_{\alpha} &= \sup_{t\in[0,t_f]} \big[\|y(t)\|\exp(-\a (t_f-t))\big], & \text{for}\ \ y\in\mathcal{Y}, \\
	\| z \|_{\a} &= \sup_{t\in[0,t_f]}\big[\|z(t)\|\exp(-\a (t_f-t))\big], & \text{for}\ \ z\in\mathcal{Z}, 
\end{align*}
in which $\|v\|=\sum_{i=1}^n |v_i|$ for any $v\in\mathbb{R}^n$ and $\|D\| = \max_{1\leq j\leq n} \sum_{i=1}^n |D_{ij}|$ for $D\in \mathbb{R}^{n\times n}$, and the parameter $\alpha$ 
serves as an additional degree of freedom to control the rate of convergence \cite{Hofer1988}, we define the operators $T_1:\mathcal{X} \times \mathcal{Y} \times \mathcal{Z} \rightarrow \mathcal{X}$, $T_2: \mathcal{X} \times \mathcal{Y} \times \mathcal{Z} \rightarrow \mathcal{Y}$, and $T_3:\mathcal{X} \times \mathcal{Y} \times \mathcal{Z} \rightarrow \mathcal{Z}$ that describe the dynamics of $x\in\mathcal{X}$, $K\in\mathcal{Y}$, and $s\in\mathcal{Z}$ as described in \eqref{eq:x_k}, \eqref{eq:K}, and \eqref{eq:Sk} by
\begin{align}
	\label{eq:T1}
	\frac{d}{dt} & T_1[x,K,s](t) = \tilde{A}(x(t),K(t),s(t))T_1[x,K,s](t) + g\\
	& -\tilde{O}(x(t))\Big(T_3[x,K,s](t)+T_2[x,K,s](t)T_1[x,K,s](t)\Big),  \nonumber \\
	&T_1[x,K,s](0) = x_0, \nonumber \\
	\frac{d}{dt}&T_2[x,K,s](t) = - T_2[x,K,s](t)\tilde{A}(x(t),K(t),s(t))\nonumber \\
	\label{eq:T2}
	& - \tilde{A}^T(x(t),K(t),s(t))T_2[x,K,s](t) \\
	& + T_2[x,K,s](t)\tilde{O}(x(t))T_2[x,K,s](t), \nonumber  \\
	&T_2[x,K,s](t_f) = 2I_n, \nonumber \\
	\label{eq:T3}
	\frac{d}{dt}&T_3[x,K,s](t) = -\tilde{A}^T(x(t),K(t),s(t))T_3[x,K,s](t) \\
	& + T_2[x,K,s](t)\tilde{O}(x(t))T_3[x,K,s](t) -  T_2[x,K,s](t) g, \nonumber  \\
	& T_3[x,K,s](t_f) = -2x_d. \nonumber
\end{align}
With these definitions, the convergence of the iterative method can be illustrated using the fixed-point theorem.

\begin{theorem}
	\label{thm:convergence}
	Consider the iterative method applied to the optimal control problem \eqref{eq:oc2} with the iterations evolving according to 
	\begin{align}
		\label{eq:T1k}
		x^{(k+1)}(t) &= T_1[x^{(k)},K^{(k)},s^{(k)}](t),  \\
		\label{eq:T2k}
		K^{(k+1)}(t) &= T_2[x^{(k)},K^{(k)},s^{(k)}](t),  \\
		\label{eq:T3k}
		s^{(k+1)}(t) &= T_3[x^{(k)},K^{(k)},s^{(k)}](t),
	\end{align}
	where the operators 
	$T_1$, $T_2$, and $T_3$ are defined in \eqref{eq:T1}, \eqref{eq:T2}, and \eqref{eq:T3}, respectively. These mappings are contractive, namely, there exists some matrix $M\in\mathbb{R}^{3\times 3}$ with all eigenvalues within a unit sphere, such that for all $x^{(k)}\in C([0,t_f];\, \mathbb{R}^n)$, $K^{(k)}\in C([0,t_f];\, \mathbb{R}^{n\times n})$, and $s^{(k)}\in C([0,t_f];\, \mathbb{R}^n)$, $k=1,2,\ldots$, the inequality
	\begin{align}
		\label{eq:contraction}
		& \left[\begin{array}{c}\| T_1[x^{(k)},K^{(k)},s^{(k)}]-T_1[x^{(k-1)},K^{(k-1)},s^{(k-1)}]\|_{\a} \\
		\| T_2[x^{(k)},K^{(k)},s^{(k)}]-T_2[x^{(k-1)},K^{(k-1)},S^{(k-1)}]\|_{\a} \\
		\|T_3[x^{(k)},K^{(k)},s^{(k)}]-T_3[x^{(k-1)},K^{(k-1)},s^{(k-1)}]\|_{\a} \end{array}\right] \nonumber \\
		& \leq M \left[\begin{array}{c}\|x^{(k)}-x^{(k-1)}\|_{\a} \\
		\|K^{(k)}-K^{(k-1)}\|_{\a} \\
		\|s^{(k)}-s^{(k-1)}\|_{\a} \end{array}\right]
	\end{align}
	holds component-wise, if 
	the magnitude of $R$ is chosen sufficiently large such that 
	\begin{equation} 
		\label{eq:criterion}
		|m_{11}|+|m_{21}|+|m_{31}|<1,
	\end{equation}
	where $m_{ij}$ is the $ij^{th}$ entry of $M$. Consequently, starting with a triple of feasible trajectories $(x^{(0)},K^{(0)},s^{(0)})$, the iteration procedure is convergent with the sequences $\{x^{(k)}\}$, $\{K^{(k)}\}$, and $\{s^{(k)}\}$ converging to the respective fixed points, $x^*$, $K^*$, and $s^*$, i.e., $x^{(k)}\to x^*$, $K^{(k)}\to K^*$, and $s^{(k)}\to s^*$.
\end{theorem}

{\it Proof:} 
From \eqref{eq:Ak} and \eqref{eq:Bk}, for each fixed $t\in [0,t_f]$, we obtain the bounds,
\begin{align*}
& \|\Delta A^{(k+1)}\| \leq [\sum_{i=1}^n \|P_i \| ^2 ] ^{1/2} \| p^{(k+1)}-p^{(k)}\| + [\sum_{i,j=1}^n \|Q_{ij} \| ^2 ] ^{1/2} \\
& \qquad \cdot ( \| p^{(k+1)} \| \| x^{(k+1)}-x^{(k)} \| + \| x^{(k)} \| \| p^{(k+1)}-p^{(k)}\| ), \\ 
& \|\Delta O^{(k+1)}\| \leq [\sum_{i,j=1}^n \|Q_{ij} \| ^2 ] ^{1/2} \| (x^{(k+1)})^2 - (x^{(k)})^2 \|,
\end{align*}
where $\Delta A^{(k+1)}\doteq\tilde{A}^{(k+1)}-\tilde{A}^{(k)}$, 
$\Delta O^{(k+1)}\doteq\tilde{O}^{(k+1)}-\tilde{O}^{(k)}$, $P_i = N_i R^{-1} B^T + B R^{-1} N^T_i$, and $Q_{ij}=N_i R^{-1} N^T_j + N_j R^{-1} N^T_i$ for $i,j = 1, \cdots, n$. We can further write these bounds in terms of $\Delta x^{(k+1)}\doteq x^{(k+1)}-x^{(k)}$, $\Delta K^{(k+1)}\doteq K^{(k+1)}-K^{(k)}$, and $\Delta s^{(k+1)}\doteq s^{(k+1)}-s^{(k)}$, given by
\begin{align} 
	& \|\Delta A^{(k+1)}\| \leq \big\{ \Big[\sum_{i=1}^n \|P_i \| ^2 \Big] ^{1/2} + \|x^{(k)} \|\Big[\sum_{i,j=1}^n \|Q_{ij} \| ^2 \Big] ^{1/2} \big\} \, \cdot \nonumber \\
	\label{eq:dA}
	&\ \big\{ \| K^{(k+1)} \|\|\Delta x^{(k+1)}\| + \|\Delta K^{(k+1)}\| \|x^{(k)} \| + \|\Delta s^{(k+1)}\| \big\} \\
	&\  + \big[\sum_{i,j=1}^n \|Q_{ij} \| ^2 \big] ^{1/2} ( \| K^{(k+1)} \| \| x^{(k+1)} \| + \| s^{(k+1)} \| ) \|\Delta x^{(k+1)}\|, \nonumber \\
	\label{eq:dO}
	& \|\Delta O^{(k+1)}\|\leq\Big[\sum_{i,j=1}^n \|Q_{ij} \|^2 \Big]^{1/2} \Big\{\| x^{(k+1)} \| + \| x^{(k)} \| \Big\} \|\Delta x^{(k+1)}\|. 
\end{align}
In addition, using \eqref{eq:Sk} we obtain
\begin{align*}
\dfrac{d}{dt}\Delta s^{(k+1)}& =  -\big[\tilde{A}^{(k)} - \tilde{O}^{(k)}K^{(k+1)} \big]^T \Delta s^{(k+1)} - \Delta K^{(k+1)}g \\
& +\big\{ \tilde{O}^{(k-1)} \Delta K^{(k+1)}+ \Delta O^{(k)} K^{(k+1)}-\Delta A^{(k)} \big\} ^T s^{(k)} 
\end{align*}
with $\Delta s^{(k+1)}(t_f) = 0$.
Then, we have
\begin{align}
	& \Delta s^{(k+1)}(t) \leq \int_t^{t_f} \Big[\beta_1 \, \| \Delta A^{(k)}(\s)\| + \b_2 \, \|\Delta O^{(k)}(\s)\| \nonumber\\
	\label{eq:dS}	& + \b_3 \, \| \Delta K^{(k+1)}(\sigma) \| \Big] d\sigma{\color{blue},}
\end{align}
where $\beta_1$, $\beta_2$ and $\beta_3$ are finite time-varying coefficients (see Appendix \ref{appd:betas}). Using \eqref{eq:K}, we can write the differential equation for $\Delta K^{(k+1)}$, that is, 
\begin{align}
	\dfrac{d}{dt} \Delta K^{(k+1)}&  =- \Delta K^{(k+1)} \Big[\tilde{A}^{(k)} - \tilde{O}^{(k)} K^{(k+1)} \Big] \nonumber\\
	\label{eq:Kd}
	& - \Big[\tilde{A}^{(k-1)} - \tilde{O}^{(k-1)}K^{(k)} \Big]^T \Delta K^{(k+1)} \\
	& - K^{(k)}\Delta A^{(k)} - (\Delta A^{(k)} )^T K^{(k+1)} + K^{(k)}\Delta O^{(k)}  K^{(k+1)}, \nonumber
\end{align}
with the terminal condition $\Delta K^{(k+1)}(t_f)=0$. Applying the variation of constants formula to \eqref{eq:Kd}, backward in time from $t=t_f$, and evaluating the norm yield the inequality
\begin{align} 
	\label{eq:dK}
	& \|\Delta K^{(k+1)}(t)\| \leq \int_t^{t_f} \Big[\beta_4 \| \Delta A^{(k)}(\s)\|+\beta_5 \| \Delta O^{(k)}(\s)\| \Big] d\s,
\end{align}
where $\beta_4$ and $\beta_5$ are both finite time-varying coefficients (see Appendix \ref{appd:betas}). 

Similarly, from \eqref{eq:x_k} and \eqref{eq:p_k}, we can write the differential equation for $\Delta x^{(k+1)}$, that is,
\begin{align*}
\frac{d}{dt} \Delta x^{(k+1)} & = \big[\tilde{A}^{(k)} - \tilde{O}^{(k)}K^{(k+1)} \big] \Delta x^{(k+1)} \\
& + \big\{ \Delta A^{(k)} - \Delta O^{(k)} K^{(k+1)} - \tilde{O}^{(k-1)} \Delta K^{(k+1)}\big\} x^{(k)} \\
& -\Delta O^{(k)} s^{(k+1)}-\tilde{O}^{(k-1)}\Delta s^{(k+1)}, 
\end{align*}
with the initial condition $\Delta x^{(k+1)}(0)=0$, 
which leads to 
\begin{align}
	& \| \Delta x^{(k+1)}(t)\| \nonumber \leq \int_0^t \Big[ \beta_6\| \Delta A^{(k)}(\s)\| + \beta_7 \| \Delta K^{(k+1)}(\s)\| \nonumber \\
	\label{eq:dx}
	& + \beta_8 \| \Delta O^{(k)}(\s)\|+ \beta_9 \|\Delta s^{(k+1)}(\s)\| \Big] d\s{\color{blue},} 
\end{align}
where $\beta_6$, $\beta_7$, $\beta_8$ and $\beta_9$ are finite time-varying coefficients (see Appendix \ref{appd:betas}). Combining the bounds in \eqref{eq:dA}, \eqref{eq:dO}, \eqref{eq:dS}, \eqref{eq:dK}, and \eqref{eq:dx}, and using the definitions of the operators $T_1$, $T_2$, and $T_3$ in \eqref{eq:T1}, \eqref{eq:T2} and \eqref{eq:T3}, respectively, we reach the inequality as in \eqref{eq:contraction} that holds component-wise, in which all of the elements, $m_{ij}$, of $M$ depend on $R^{-1}$ through $\delta$ and $\zeta$ (see Appendix \ref{appd:M}). 

Now, observe that $m_{i1}>m_{i2}$ and $m_{i1}>m_{i3}$ for $i=1,2,3$, for example, for $i=3$, $m_{31}-m_{32}> \zeta (\beta_2 + \beta_3 \beta_5) ( \|x^{(k)}\|+\|x^{(k-1)}\| )>0$ and $m_{31}-m_{33}> \zeta (\beta_2 + \beta_3 \beta_5) ( \|x^{(k)}\|+\|x^{(k-1)}\| )>0$. Therefore, due to the multiplicative influence of $R^{-1}$, the condition \eqref{eq:criterion} can be fulfilled by choosing a sufficiently large $R$, 
which leads to $\rho(M)\leq \| M \|_1<1$, where $\rho(M)$ is the spectral radius of $M$. Now, by proper construction of the invariant sets,
\begin{align*}
D_1 & \doteq \{x\in \mathcal{X}:  \|x - x^{(0)} \|_\a \leq l_1\}, \\
D_2 & \doteq \{K\in \mathcal{Y}:  \|K - K^{(0)} \|_\a \leq l_2\}, \\
D_3 & \doteq \{s\in \mathcal{Z}:  \|s - s^{(0)} \|_\a \leq l_3\}
\end{align*}
with $l_1, l_2, l_3>0$ and appropriate choice of $R$, the contraction property of $T_1$, $T_2$, and $T_3$ in $D_1\times D_2\times D_3$ can be guaranteed. It follows that the iteration procedure is convergent by the fixed-point theorem \cite{Banach22} with the convergent fixed points $x^{(k)}\rightarrow x^*$, $K^{(k)}\rightarrow K^*$, and $s^{(k)}\rightarrow s^*$.

This iterative method can be directly applied to solve the same class of optimal control problems involving a bilinear ensemble system, which will be discussed in Section \ref{sec:ensemble}.

\begin{remark}
	\rm Note that the control weight matrix $R$ has to be large enough to compensate for the translation drifted by $g$, interpreted as the condition \eqref{eq:criterion}, so as to guarantee the convergence. This penalizes the cost of control energy, which in turn results in greater terminal cost and less control energy.
\end{remark}

\subsection{Global Optimality of the Convergent Solution}
\label{sec:optimality}
While the iterative procedure is convergent, the convergent control $u^*$ satisfying \eqref{eq:x_k} and \eqref{eq:p_k} as $k\to\infty$ will in turn satisfy \eqref{eq:state1} and \eqref{eq:costate1}, which are the necessary optimality condition for Problem \eqref{eq:oc1}. In this section, we will further illustrate that $u^*$ is a global optimal control given appropriate conditions on the value function associated with Problem \eqref{eq:oc1}.


For each iteration $k$, \eqref{eq:oc2} is a time-dependent problem with the time horizon specified, and the optimal control satisfies the Hamilton-Jacobi-Bellman (HJB) equation, given by
\begin{align} 
	\label{eq:HJBlqr}
	\nonumber
	&\min_{u\in \mathcal{U}} \Big\{ \dfrac{\partial V^{(k)}}{\partial x^{(k)}}(t,x^{(k)})^T (\tilde{A}^{(k-1)}x^{(k)}+\tilde{B}^{(k-1)}u + g) +\dfrac{1}{2} u^TRu\Big\}\\
	& + \dfrac{\partial V^{(k)}}{\partial t}(t,x^{(k)}) \equiv 0,
\end{align}
with the boundary condition $V^{(k)}(t_f,x^{(k)})=0$, where $V^{(k)}$ is the value function and $\mathcal{U}\subseteq PC([0,t_f];\, \mathbb{R}^m)$ is the set of all piecewise continuous admissible controls. 
Since the matrix $R\in \mathbb{R}^{n\times n}$ is positive definite, the function to be minimized in \eqref{eq:HJBlqr} is strictly convex in the control variable $u$. As a result, the minimization problem in \eqref{eq:HJBlqr} has a unique solution given by the stationary point, satisfying $$(\tilde{B}^{(k-1)})^T \dfrac{\partial V^{(k)}}{\partial x^{(k)}}(t,x^{(k)}) + Ru_*^{(k)}(t)=0,$$
or, equivalently,
\begin{align}
	\label{eq:optimalu}
	u_*^{(k)}(t) = -R^{-1}(\tilde{B}^{(k-1)})^T \dfrac{\partial V^{(k)}}{\partial x^{(k)}}(t,x^{(k)}).
\end{align}
Substituting \eqref{eq:optimalu} into the HJB equation in \eqref{eq:HJBlqr} gives a first-order nonlinear partial differential equation, 
\begin{align}
	& \dfrac{\partial V^{(k)}}{\partial t}(t,x^{(k)}) + \dfrac{1}{2} \Big[ \dfrac{\partial V^{(k)}}{\partial x^{(k)}}(t,x^{(k)})^T\tilde{A}^{(k-1)}x^{(k)} \nonumber \\ 
	\label{eq:fopde}
	& + (x^{(k)})^T(\tilde{A}^{(k-1)})^T\dfrac{\partial V^{(k)}}{\partial x^{(k)}}(t,x^{(k)})\Big] + \dfrac{\partial V^{(k)}}{\partial x^{(k)}}(t,x^{(k)})^T g \\
	& -\dfrac{1}{2} \dfrac{\partial V^{(k)}}{\partial x^{(k)}}(t,x^{(k)})^T \tilde{O}^{(k-1)} \dfrac{\partial V^{(k)}}{\partial x^{(k)}}(t,x^{(k)})\equiv 0. \nonumber
\end{align}
Due to the quadratic and symmetric nature of \eqref{eq:fopde}, we consider the value function of the form,
\begin{align}
	\label{eq:Vfunc}
	V^{(k)}(t,x) = \dfrac{1}{2}x^TK^{(k)}x + x^Ts^{(k)} + \dfrac{1}{2} q^{(k)},
\end{align}
with the boundary condition $V^{(k)}(t_f,x^{(k)})=0$, where $K^{(k)}(t)\in \mathbb{R}^{n\times n},~s^{(k)}(t)\in \mathbb{R}^{n}, ~\forall t\in [0,t_f]$ and  $q^{(k)}\in\mathbb{R}$ satisfies the differential equation 
$$\dot{q}^{(k)} = (s^{(k)})^T \tilde{O}^{(k)} s^{(k)} - 2(s^{(k)})^T g,$$
with $q^{(k)}(t_f) = 0$. It is straightforward to verify that $(V^{(k)}(t,x^{(k)}),u_*^{(k)})$ is a classical solution to the HJB equation \eqref{eq:fopde} if 
$K^{(k)}(t)$ and $s^{(k)}(t)$ satisfy the 
differential equations \eqref{eq:K} and \eqref{eq:Sk}, respectively, with the respective boundary conditions $K^{(k)}(t_f)=2I_n$ and $s^{(k)}(t_f)=-2x_d$. Moreover, the minimizing control for Problem \eqref{eq:oc2} is given by the linear feedback law as in \eqref{eq:optimalu}.

We now have shown that the control $u_*^{(k)}$ presented in \eqref{eq:optimalu} is the global optimum for the $k^{th}$ iteration in Problem \eqref{eq:oc2}. Next, we will show that the convergent solution $u^*$ of $u_*^{(k)}$, i.e., $u_*^{(k)}\rightarrow u^*$, for the linear problem \eqref{eq:oc2} is indeed a global optimal control for the original bilinear problem \eqref{eq:oc1} under some regularity conditions on the value function associated with \eqref{eq:oc2} presented in \eqref{eq:Vfunc}.

\begin{theorem}
	\label{thm:globalsufficient}
	Let $u^*$ be the convergent solution of the optimal control sequence $\{u_*^{(k)}\}$ generated by the iterative procedure,  
	i.e., $u_*^{(k)}\rightarrow u^*$ as $k\to\infty$, with respect to a given initial condition $x_0$, and let $V^*$ be the corresponding convergent value function defined in \eqref{eq:Vfunc}, i.e., $V^{(k)} \rightarrow V^*$. If (i) $V^*\in C^1$, and $\dfrac{\partial V^*}{\partial t}$ and $\dfrac{\partial V^*}{\partial x}$ are Lipschitz continuous; and (ii) there exist real-valued $L_1$ functions, $\xi \in L_1([0,t_f])$ and $\eta_i \in L_1(\mathbb{R}^n)$, $i = 1,2,\ldots,n$, such that $\Big|\dfrac{\partial V^{(k)}}{\partial t}(t,x^{(k)}(t))\Big| \leq \xi(t)$ for all $k\in\mathbb{N}$ and for all $t\in [0,t_f]$, and for each component $i$, $\Big| \Big[\dfrac{\partial V^{(k)}}{\partial x}(t,x)\Big]_i\Big| \leq \eta_i(x)$ for all $k\in\mathbb{N}$ and for all $x=x^{(k)}(t)\in\mathbb{R}^n$, then $u^*$ is a global optimum for the original Problem \eqref{eq:oc1}.
\end{theorem}

{\it Proof:} See Appendix \ref{appd:proof}.

\section{Optimal Control of Inhomogeneous Bilinear Ensembles}
\label{sec:ensemble}
In this section, we extend the iterative method presented in Section \ref{sec:inhomogeneous} to deal with optimal control problems involving a 
bilinear ensemble system as in \eqref{eq:bilinear_ensemble}. 
Specifically, we study the ensemble analogy of the optimal control problem \eqref{eq:oc1}, given by
\begin{align}	
	\label{eq:oc3}
	\min & \ \  \mathcal{J} = \frac{1}{2}\int_0^{t_f} u^T(t)Ru(t) \, dt + J_1,  \tag{P3}\\ 
	{\rm s.t.} & \ \ \frac{d}{dt}{X(t,\b)}=A(\b)X+B(\b)u+\big[\sum_{i=1}^m u_i B_i(\b)\big] X + g(\b), \nonumber \\
	& \quad \ \ X(0,\b) = X_0(\b), \nonumber
\end{align}
where $X_0\in C(K;\mathbb{R}^n)$ is a continuous function over $K$, the bilinear system $X(t,\b)$ is defined similarly as the system in \eqref{eq:bilinear_ensemble} with $g$ continuous over $K$, the control weight matrix $R\in\mathbb{R}^{m\times m}$ is positive definite, and the terminal cost $J_1$ is defined as the averaged terminal error over $\b$ by 
\begin{align}
	\label{eq:J_1}
	J_1 = \dfrac{1}{\mu (K)} \int_K \|X(t_f,\b)-X_d(\b)\|_K^2 \, d\mu(\b),
\end{align}
where $X_d\in C(K;\mathbb{R}^n)$ and $\mu(K)$ is the Lebesgue measure of $K$.

Following the iterative procedure developed in Section \ref{sec:iterative}, we represent the bilinear ensemble system in \eqref{eq:oc3} as an iteration equation and formulate the optimal ensemble control problem as
\begin{align}	
	\label{eq:oc4}
	\min & \mathcal{J} = \frac{1}{2}\int_0^{t_f} u^T(t)Ru(t) \, dt + J_1,  \tag{P4}\\ 
	{\rm s.t.} & \frac{d}{dt}{X^{(k)}(t,\b)}=\tilde{A}^{(k-1)}(t,\b)X^{(k)}+\tilde{B}^{(k-1)}(t,\b) u^{(k)}+g(\b) \nonumber \\
	& \ \ \ X^{(k)}(0,\b)=X_0(\b), \nonumber
\end{align}
where the matrices $\tilde{A}^{(k-1)}(t,\b)\in\mathbb{R}^{n\times n}$ and $\tilde{B}^{(k-1)}(t,\b)\in\mathbb{R}^{n\times m}$ for $t\in[0,t_f]$ and $\b\in K$ are defined analogously as in \eqref{eq:Ak} and \eqref{eq:Bk}.

\subsection{Iterative Method for Optimal Ensemble Control}
Now, we validate the applicability of the iterative method presented in Sections \ref{sec:iterative} and \ref{sec:convergence} for solving the optimal ensemble control problem \eqref{eq:oc3}. 
In particular, we will demonstrate that following the iterative method the optimal solution of the discretized problem of \eqref{eq:oc4} (discretization over the parameter domain $K$), which is in the form of \eqref{eq:oc2}, converges to that of the discretized problem of \eqref{eq:oc3} (provided the existence of the optimal solution), which is in the form of \eqref{eq:oc1}.

 \begin{figure}[t]
  	\centering
  	\includegraphics[width=0.9\columnwidth]{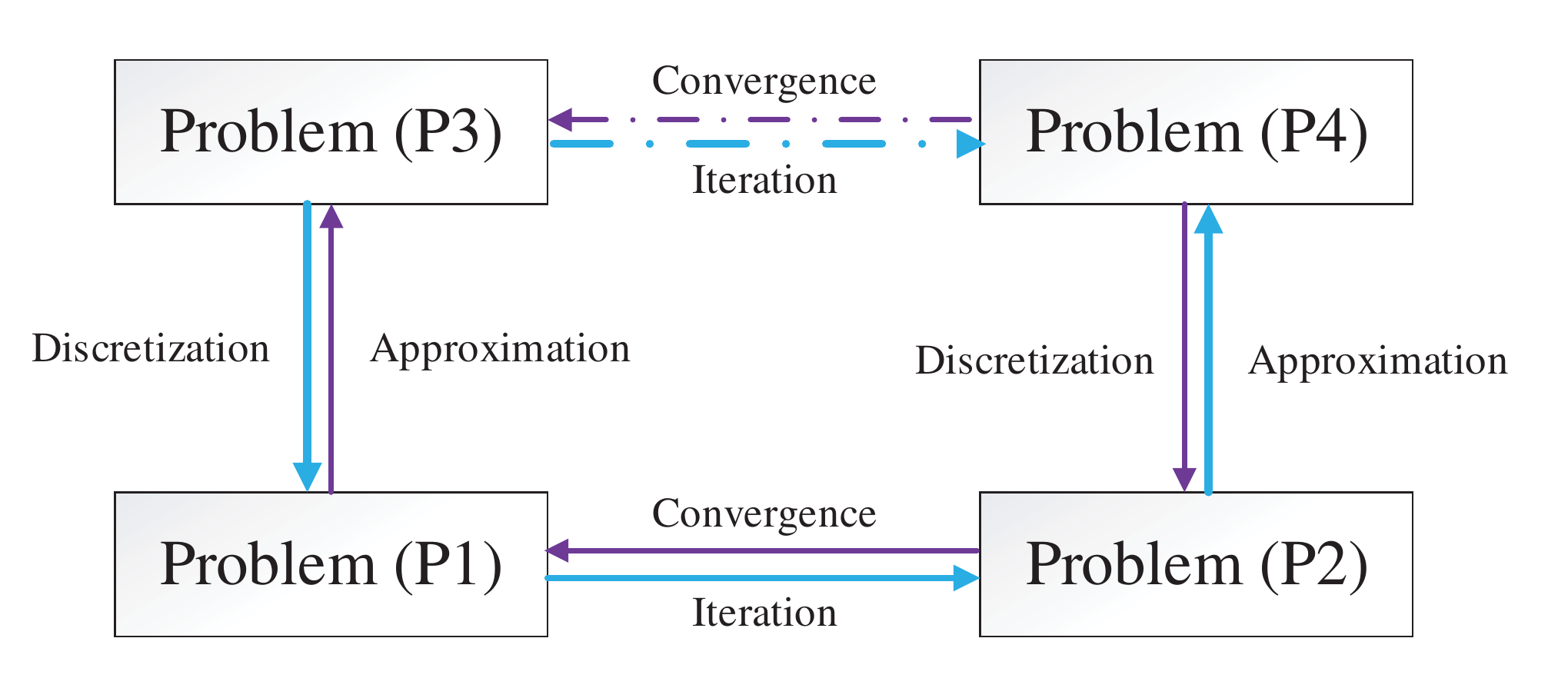}
  	\caption{A schematic diagram for the relationship between the four types of optimal control problems.} 
 	\label{fig:Diagram}
  \end{figure}

Our validation is based on the following steps:

\eqref{eq:oc3}$\Rightarrow$\eqref{eq:oc1} (Discretization): We first consider a finite uniform sample, $\{\b_1,\ldots,\b_q\}\subset K$, $q\in\mathbb{Z}^+$, of the ensemble $X(t,\b)$ in \eqref{eq:oc3} with $\b_1=a$ and $\b_q=b$, and then formulate Problem \eqref{eq:oc1} involving the state vector $\hat{X}(t)=[X(t,\b_1),\ldots,X(t,\b_q)]^T\in\mathbb{R}^{nq}$, $A=A_1\oplus\cdots\oplus A_q$, and $B=B_1\oplus\cdots\oplus B_q$, where $A_i=A(\b_i)$ and $B_i=B(\b_i)$ for $i=1,\ldots,q$. Then, the cost functional is $J = \frac{1}{2}\int_0^{t_f} u^T(t)Ru(t) \, dt + \frac{1}{q}\|\hat{X}(t_f) - \hat{X}_d\|_2^2$, where $\hat{X}_d = [X_d(\b_1),\ldots, X_d(\b_q)]^T$ denotes the target state.

\eqref{eq:oc1}$\Rightarrow$\eqref{eq:oc2} (Iterations): Then, following the iterative method developed in Sections \ref{sec:iterative} and \ref{sec:convergence}, we can solve \eqref{eq:oc1} by solving \eqref{eq:oc2}, and find, at each iteration $k$, the optimal solution $\hat{X}^{(k)}$ to the problem as in \eqref{eq:oc2}, given by
\begin{align}
	\label{eq:P2}
	\min &\quad J_q=\frac{1}{2}\int_0^{t_f} (u^{(k)})^T(t)Ru^{(k)}(t) \, dt + \dfrac{1}{q}\Vert \hat{X}^{(k)}(t_f) - \hat{X}_d \Vert _2^2, \nonumber\\
	{\rm s.t.} &\quad \dot{\hat{X}}^{(k)} =\tilde{A}^{(k-1)} \hat{X}^{(k)}+\tilde{B}^{(k-1)} u^{(k)} + g, \\  
	&\quad \hat{X}^{(k)}(0)=\hat{X}_0,  \nonumber
\end{align}
where $\hat{X}^{(k)}(t)=[X^{(k)}(t,\b_1),\ldots,X^{(k)}(t,\b_q)]^T\in\mathbb{R}^{nq}$ and $\hat{X}_0 = [X_0(\b_1),\ldots, X_0(\b_q)]^T$ 
denotes the initial state. 

\eqref{eq:oc2}$\Rightarrow$\eqref{eq:oc4} (Approximation): Taking the Riemann sum of the terminal error, we have, by the continuity of $X^{(k)}(t_f,\b)$,
\begin{align}
	\label{eq:Euler}
	& \dfrac{1}{q}\Vert \hat{X}^{(k)}(t_f) - \hat{X}_d \Vert _2^2 = \dfrac{1}{q} \sum_{j=1}^q \Vert X^{(k)}(t_f,\b_j) - X_d(\b_j) \Vert _2^2 \nonumber \\
	& \rightarrow \dfrac{1}{\mu (K)} \int_K \|X^{(k)}(t_f,\b)-X_d(\b)\|^2 \, d\mu(\b) = J_1
\end{align} 
as $q\rightarrow\infty$. This implies that, for each $k$, the solution to the optimal control problem in \eqref{eq:P2} (as in the form of \eqref{eq:oc2}) approaches that of Problem \eqref{eq:oc4} as $q\to\infty$.

On the other hand, now we show that the reversed direction, i.e., \eqref{eq:oc4}$\Rightarrow$\eqref{eq:oc2}$\Rightarrow$\eqref{eq:oc1}$\Rightarrow$\eqref{eq:oc3} in Figure \ref{fig:Diagram}, holds by the convergence of the iterative algorithm, and then conclude the extension of the iterative method for solving optimal ensemble control problem \eqref{eq:oc3} via solving \eqref{eq:oc4} iteratively. Note it is in practice via solving the discretized problem of \eqref{eq:oc4}).

\eqref{eq:oc4}$\Rightarrow$\eqref{eq:oc2} (Discretization): Starting with Problem \eqref{eq:oc4}, we consider a discretized problem as in \eqref{eq:P2} (in the form of \eqref{eq:oc2}) through uniform discretization over $K$, in which, at each iteration $k$, the state is $\hat{X}^{(k)}(t)=[X^{(k)}(t,\b_1),\ldots,X^{(k)}(t,\b_q)]^T\in\mathbb{R}^{nq}$ and the system matrices are $\tilde{A}^{(k)}=\tilde{A}_1^{(k)}\oplus\cdots\oplus\tilde{A}_q^{(k)}\in\mathbb{R}^{nq\times nq}$ and $\tilde{B}^{(k)}=\tilde{B}_1^{(k)}\oplus\cdots\oplus\tilde{B}_q^{(k)}\in\mathbb{R}^{nq\times mq}$ with $\tilde{A}_i^{(k)}=\tilde{A}^{(k)}(t,\b_i)$ and $\tilde{B}_i^{(k)}=\tilde{B}^{(k)}(t,\b_i)$ for $i=1,\ldots,q$.

\eqref{eq:oc2}$\Rightarrow$\eqref{eq:oc1} (Convergence): Given the convergence of the iterative method, the optimal solution $\hat{X}^{(k)}$ of \eqref{eq:P2} (in the form of \eqref{eq:oc2}) converges to the optimal solution, $\hat{X}^*$, of the bilinear problem in the form of \eqref{eq:oc1}, i.e., $\hat{X}^{(k)}\to\hat{X}^*$ as $k\to\infty$.

\eqref{eq:oc1}$\Rightarrow$\eqref{eq:oc3} (Approximation): 
The Lipschitz continuity of $A$ and $B$ in $K$ guarantees that the terminal state $X(t_f,\b)$ is piecewise continuous in $K$, and so is the optimal solution $X^*(t_f,\b)$ to Problem \eqref{eq:oc3} provided it exits. Then, a similar argument as in \eqref{eq:Euler} can be made to claim the approximation of the solution of \eqref{eq:oc1} to that of \eqref{eq:oc3}, which follows
\begin{align*}
	& \dfrac{1}{q}\Vert \hat{X}^*(t_f) - \hat{X}_d \Vert _2^2 = \dfrac{1}{q} \sum_{j=1}^q \Vert X^*(t_f,\b_j) - X_d(\b_j) \Vert _2^2 \nonumber \\
	& \rightarrow \dfrac{1}{\mu (K)} \int_K \|X^*(t_f,\b)-X_d(\b)\|^2 \, d\mu(\b) = J_1.
\end{align*}

\begin{remark} (Application of the Iterative Method to Time-Varying Bilinear Systems)
	\rm The iterative method can be directly applied to solve for the same class of quadratic optimal control problems involving a time-varying bilinear ensemble system. The convergence and optimality of the convergent solution follow the derivations presented in Sections \ref{sec:inhomogeneous} and \ref{sec:ensemble}.
\end{remark}

\subsection{Optimal Control of Stochastic Bilinear Ensembles}
\label{sec:stochastic}
In this section, we illustrate a straightforward application of the iterative method to find optimal controls for stochastic bilinear ensemble systems in the presence of additive noise.

Consider an ensemble of independent, structurally identical stochastic bilinear systems driven by an additive noise process in the It\^{o}'s form, 
\begin{align} 
	\label{eq:SBLS}
	dX(t,\b) & =  A(\b)X(t,\b)dt + B(\b)u(t)dt \\
	& +  \Big[\sum_{i=1}^m u_i(t)B_i(\beta) \Big]X(t,\beta)dt + G(\b)dS(t), \nonumber
\end{align}
where $G\in C(K;\mathbb{R}^{n\times k})$ and the other terms in system \eqref{eq:SBLS} are defined as in \eqref{eq:bilinear_ensemble}. The term $dS$ is the differential of a continuous-time stochastic process $S(t)\in\mathbb{R}^k$. If the system $X(t,\b)$ is driven by the Poisson jump processes, i.e., $S(t)=N(t)\in\mathbb{R}^k$, where $N(t)$ consists of independent and identically distributed (i.i.d.) Poisson counters with constant intensities $\lambda\in\mathbb{R}^k$ satisfying $\E N(t)=\lambda t$ and $\E[N(t)N^T(s)]=\Lambda\min(t,s)$ with $\Lambda=\text{diag}(\lambda)$, then the evolution of the expected value of the state follows the dynamic equation
\begin{align}
	\label{eq:dE_Possion}
	\frac{d}{dt}\E X(t,\b) = & A(\b)\mathcal{E}X(t,\b)+ B(\b)u(t) \\
	&+ \Big( \sum_{i=1}^m u_i(t)B_i(\b) \Big)\mathcal{E}X(t,\b) +G(\b)\l, \nonumber
\end{align}
which is a deterministic bilinear ensemble system of the form as in \eqref{eq:oc3} with $g(\b)=G(\b)\l$. 
If $X(t,\b)$ is driven by the standard and independent Wiener processes, i.e., $S(t)=W(t)\in\mathbb{R}^k$ with $\E W(t)=0$ and $\E [W(t)W^T(s)]=I_k\min(t,s)$, where $I_k$ is the $k\times k$ identity matrix, then $g(\b)=0$.  

Therefore, the iterative method can be directly adopted to solve the stochastic optimal ensemble control problem minimizing the quadratic cost functional that considers the trade-off between the expected terminal error and the energy of the control input, given by $\mathcal{J} = \frac{1}{2}\int_0^{t_f} u^T(t)Ru(t) dt + J_E$, subject to the system dynamics described in \eqref{eq:dE_Possion} and the initial condition $X(0,\b)=X_0(\b)$, where the terminal cost $J_E$ measures the error in the mean, defined by
\begin{align}
	\label{eq:Jcost}
	J_E = \dfrac{1}{\mu (K)} \int_K \| \E X(t_f,\b)-X_d(\b)\|_K^2 \, d\mu(\b).
\end{align}

\section{Examples}
\label{sec:example}
In this section, we introduce several practical control design problems in the areas of neuroscience and quantum control as examples to illustrate the robustness and applicability of the developed iterative method. In particular, we illustrate how the choices of the matrix $R$ determine the convergence by numerical demonstration of the criterion in \eqref{eq:criterion} using a pulse design problem for coherence transfer in protein NMR spectroscopy.

\begin{figure}[t]
	\centering
  	\subfigure[The convergent optimal control]{
  	\includegraphics[width=1\columnwidth]{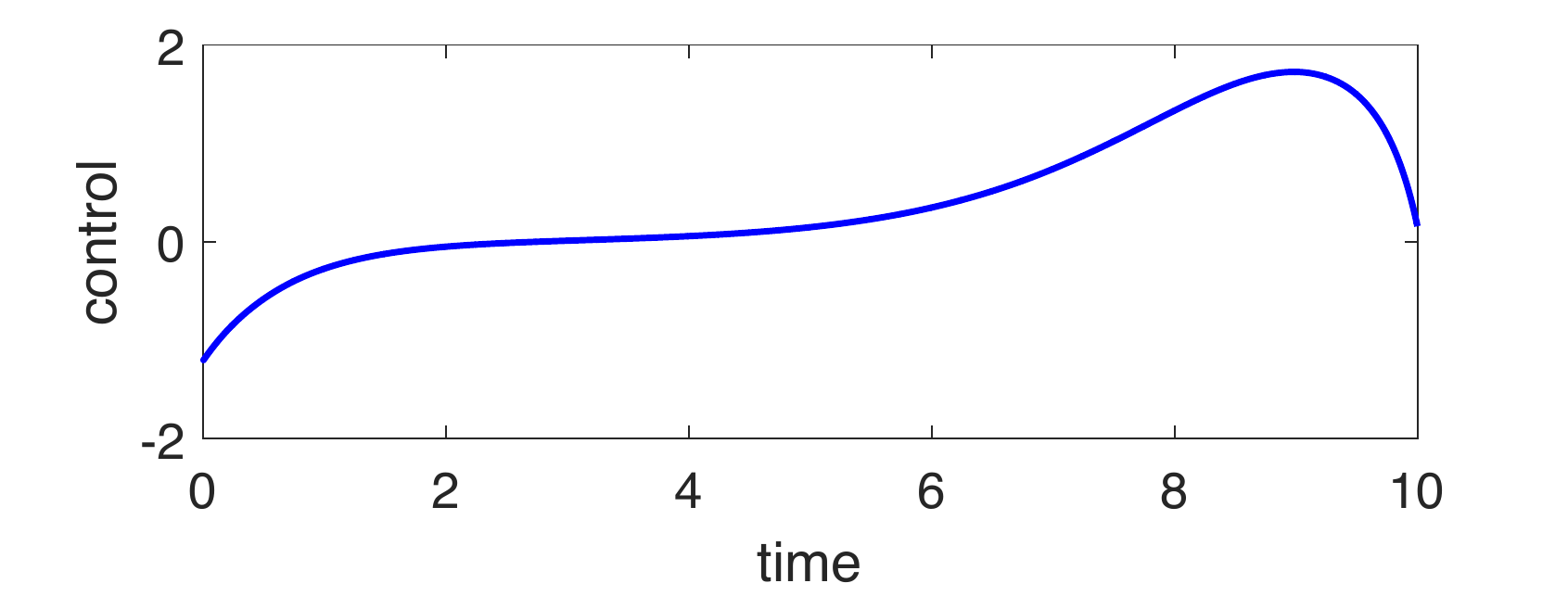} 
 	\label{fig:NControlEN}}
	\subfigure[A sample of optimal expected trajectories]
	{\includegraphics[width=1\columnwidth]{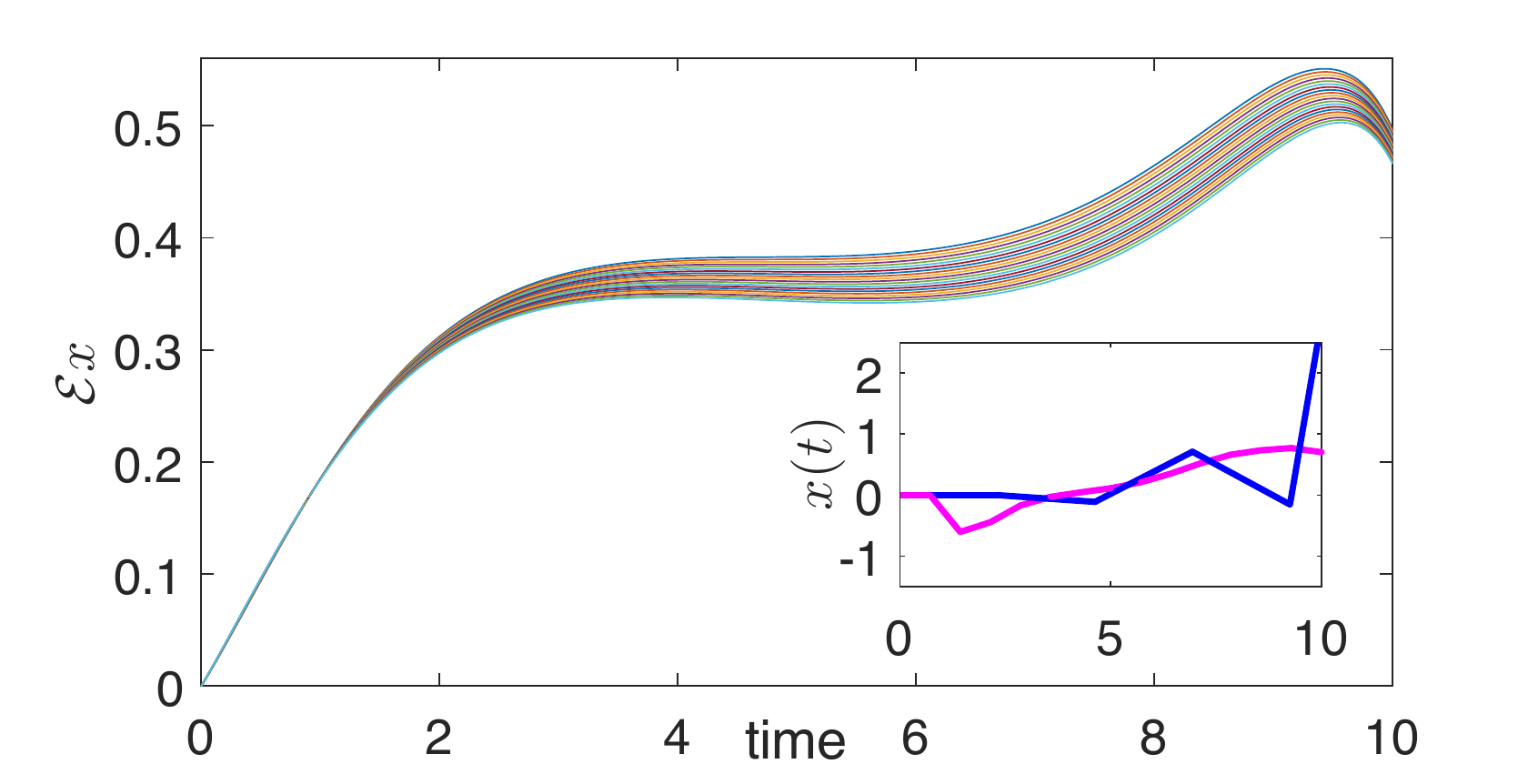}
 	\label{fig:NStateEN}}
  	\caption{\subref{fig:NControlEN}  The optimal ensemble control that spikes the IAF neuron ensemble with a variation in the decay rate $\a\in [1.2,1.3]$ modeled as in \eqref{eq:IAF} from $\mathcal{E}x(0)=0$ to $\mathcal{E}x(10)= 0.5$, while minimizing the cost functional $\mathcal{J}$ in \eqref{eq:J}. The optimal control is calculated using the iterative method that converges in 17 iterations with $R=5$ based on the stopping criterion $\|x^{(k+1)} - x^{(k)}\|< 10^{-12}$. \subref{fig:NStateEN} A sample of optimal expected trajectories following the optimal ensemble control shown in \subref{fig:NControlEN} with 2 sample paths corresponding to $\a=1.2$ and $\a=1.3$ shown in the inset.}
 	\label{fig:NeuronEnsemble}
\end{figure}

\begin{figure}[t]
	\centering
  	\subfigure[The convergent optimal control]{
  	\includegraphics[width=1\columnwidth]{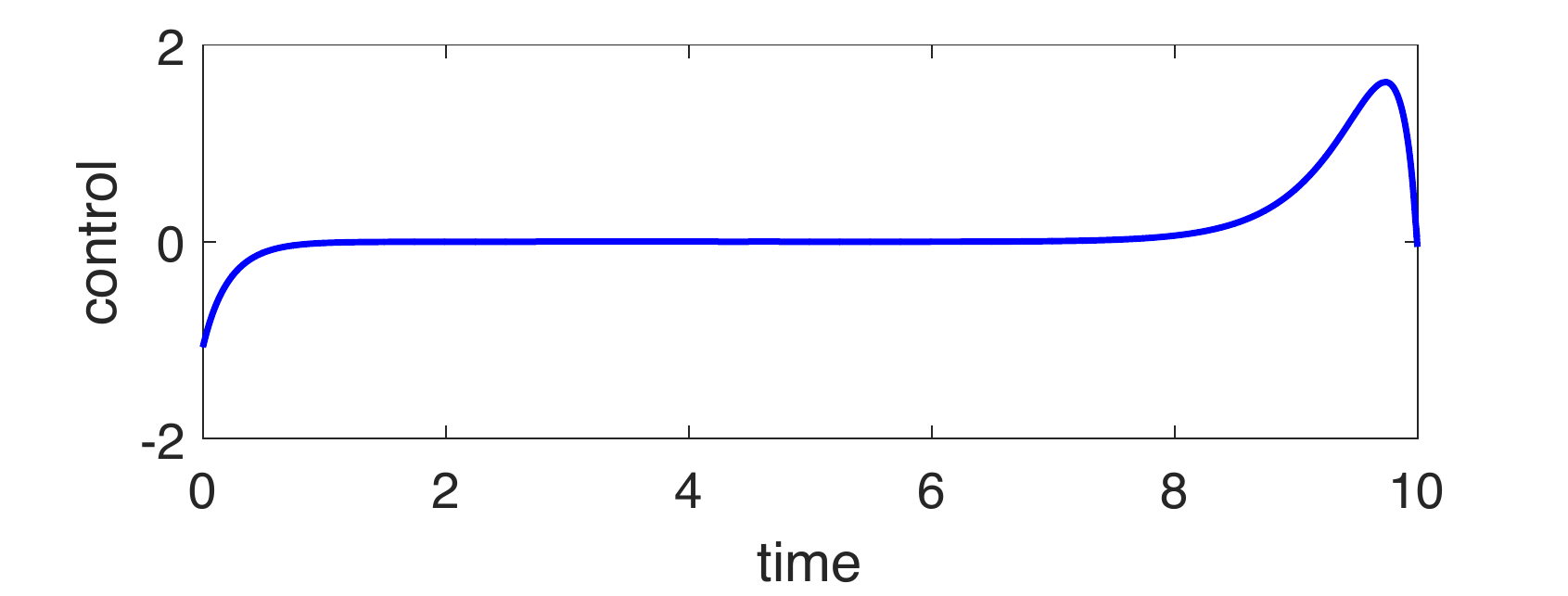} 
 	\label{fig:NControlEN2}}
	\subfigure[A sample of optimal expected trajectories]{\includegraphics[width=1\columnwidth]{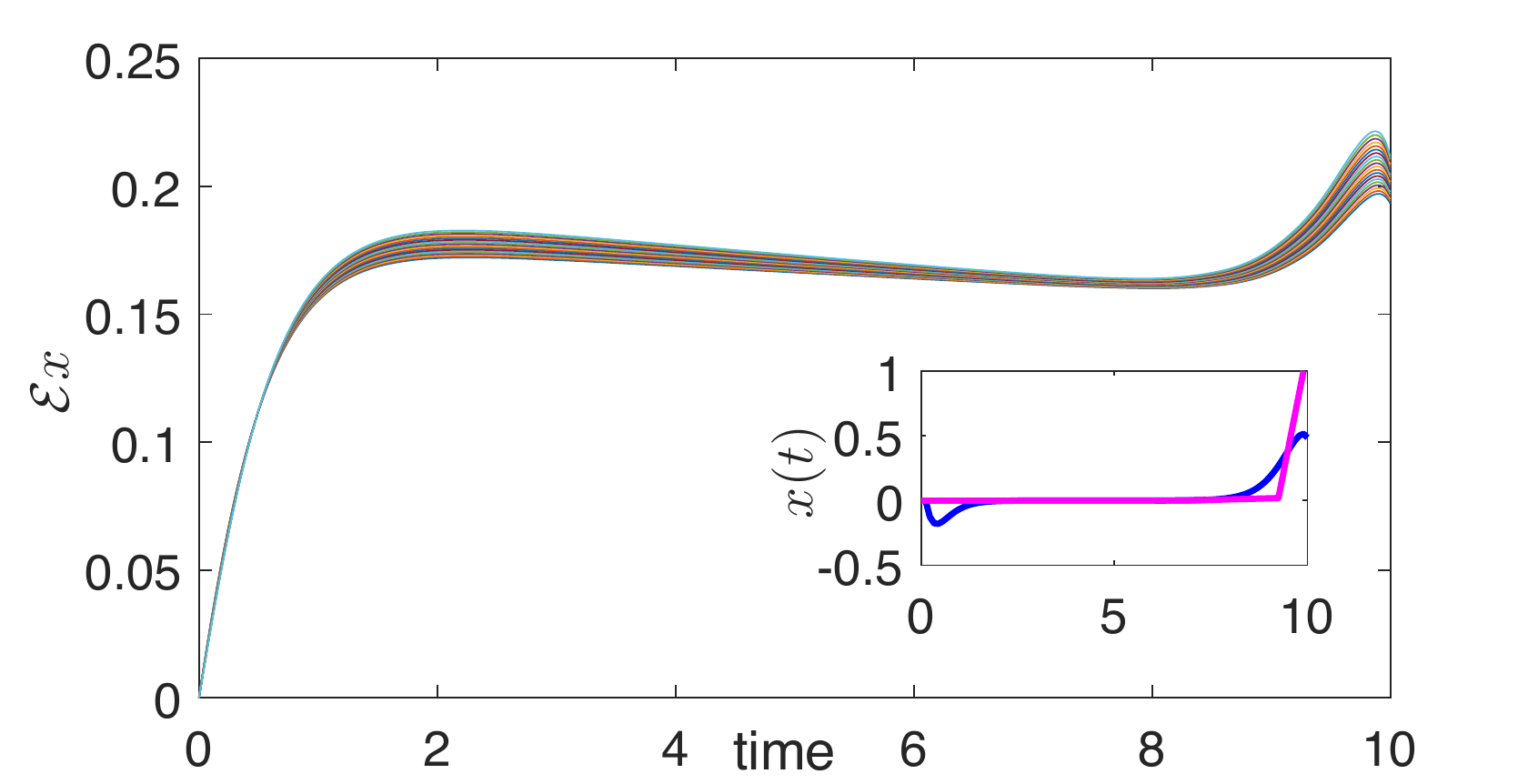}
 	\label{fig:NStateEN2}}
  	\caption{\subref{fig:NControlEN2} The optimal ensemble control that spikes the IAF neuron ensemble with a variation in the firing rate $\g\in [1.8,2.2]$ modeled as in \eqref{eq:IAF} from $\mathcal{E}x(0)=0$ to $\mathcal{E}x(10)= 0.2$, while minimizing the cost functional $\mathcal{J}$ in \eqref{eq:J}. The optimal control is calculated using the iterative method that converges in 10 iterations with $R=3$ based on the stopping criterion $\|x^{(k+1)} - x^{(k)}\| < 10^{-12}$. \subref{fig:NStateEN} A sample of optimal expected trajectories following the optimal ensemble control shown in \subref{fig:NControlEN} with 2 sample paths corresponding to $\g=1.8$ and $\g=2.2$ shown in the inset.}
 	\label{fig:NeuronEnsemble2}
\end{figure}

\begin{example}[Spiking IAF Neurons]
\label{ex:neuron}
\rm
The Integrate and Fire (IAF) neuron model is widely used in computational neuroscience to describe the dynamics of membrane potential of a population of neurons expressing Channelrhodopsin-2 (ChR2) \cite{Ching2013a}. It is of practical interest to find an optimal stimulus that controls the spikes of IAF neurons in the presence of noise. Here, we consider optimal firing of an IAF neuron population parameterized by the parameters $\a$ and $\g$, modeled by
\begin{align}
	dx(t,\a,\g) = &-\a x(t,\a,\g)dt + \g u(t)(E-x(t,\a,\g))dt \nonumber\\
	\label{eq:IAF}
	&+g(\a,\g) dN,
\end{align}
where $\a\in[\a_1,\a_2]\subset\mathbb{R}^+$ depicts the dispersion in the decay rate of the neuron, i.e., the reciprocal of the passive membrane time constant, and $\g\in[\g_1,\g_2]\subset\mathbb{R}^+$ captures the variation in responsiveness to the environment; the resetting potential $E$ is assumed to be uniform across all neurons with $E > x_d$, where $x_d$ is the firing state for ChR2; $x(t,\a,\g)\in\mathbb{R}$ represents membrane potential of the neuron expressing ChR2 and the conductance $u(t)\geq 0$ denotes the optogenetics control input; $N$ is a Poisson counter with the jumping rate $\l$ characterizing sudden changes in the membrane potential due to, for example, the opening of an ion channel; and $g(\a,\g)\in\mathbb{R}$ denotes the intensity of the noise. In particular, we optimize the objective 
\begin{equation}
	\label{eq:J}
	\mathcal{J} = \frac{1}{2}\int_0^{t_f} u^T(t)Ru(t) dt + J_E
\end{equation}
that considers the tradeoff between firing performance and the demand of input energy, where $J_E$ is defined as in \eqref{eq:Jcost}.

We examine two cases in which we consider the parameter variations in either $\a$ or $\g$:

(Case I) Suppose that, in \eqref{eq:IAF}, there exists dispersion in the decay rate $\a\in K_1=[1.2,1.3]$ and $\g=2$, $E=1$, $x_d = 0.5$, $g = 0.15$, and $\l = 2$, and consider minimizing $\mathcal{J}$ for time $t_f=10$. The optimal ensemble control derived using the iterative method is shown in Figure \ref{fig:NControlEN}. Figure \ref{fig:NStateEN} displays the resulting expected optimal trajectories, i.e., $\E x(t,\a_i,2)$, for a sample of 20 neurons with $\a_i$ uniformly sampled in $K_1$, in which the inset shows two sample paths of $x(t,\a,\g)$ for $(\a,\g)=(1.2,2)$ and $(\a,\g)=(1.3,2)$ following the optimal control presented in Figure \ref{fig:NControlEN}. The iterative method converges in $17$ iterations based on the stopping criterion $\|x^{(k+1)} - x^{(k)}\| < 10^{-12}$ when choosing $R=5$. 

(Case II) Suppose that the variation appears in the firing rate with $\g\in K_2=[1.8,2.2]$, and, in \eqref{eq:IAF}, $\a=2.6$, $E = 1$, $x_d = 0.2$, $g = 0.1$, and $\l = 4$, we obtain the optimal ensemble control in 10 iterations using $R=5$, shown in Figure \ref{fig:NControlEN2}, given the stopping criterion $\|x^{(k+1)} - x^{(k)}\| < 10^{-12}$. A sample of the resulting expected trajectories is displayed in Figure \ref{fig:NStateEN2} with two sample paths of $x(t,\a,\g)$ for $(\a,\g)=(2,1.8)$ and $(\a,\g)=(2,2.2)$ illustrated in the inset.
\end{example}

\begin{example}[Excitation of Two-Level Systems]
 \label{ex:bloch_ensemble}
 \rm Designing an optimal pulse that excites a collection of two-level systems is an essential control task that enables various applications in quantum science and technology \cite{Li_PNAS11,Glaser98}. The dynamics of a two-level system obeys the Bloch equations, which forms a bilinear control system evolving on the special orthogonal group SO(3), given by
\begin{align}
	\label{eq:blochsys}
 	\frac{d}{dt} \begin{bmatrix}
 	x_1 \\ x_2 \\ x_3
 	\end{bmatrix} = \begin{bmatrix}
 	0 & -\omega & u_1 \\
 	\omega & 0 & -u_2 \\
 	-u_1 & u_2 & 0 \end{bmatrix}
 	\begin{bmatrix}
 		x_1 \\ x_2 \\ x_3
 	\end{bmatrix},
\end{align}
where $X=(x_1,x_2,x_3)^T$ denotes the bulk magnetization of the spins, $\w$ denotes the Larmor frequency of the spins, and $u_1$ and $u_2$ are the radio-frequency (RF) fields applied on the $y$ and the $x$ direction, respectively \cite{Li_TAC09}. A typical problem in quantum control is to develop a broadcast control field, the so-called broadband pulse, driving an ensemble of systems as modeled in \eqref{eq:blochsys} with $\w\in [\w_1,\w_2]$ from the equilibrium state $X(0,\w)=X_0(\w) = (0,0,1)^T$ as close as possible to a desired excited state, e.g., $X_d(\w)=(1,0,0)^T$ at a specified time $t_f$, with minimum energy \cite{Li_PNAS11}.

\begin{figure}[t]
	\centering
	\subfigure[The optimal ensemble control]{\includegraphics[width=\columnwidth]{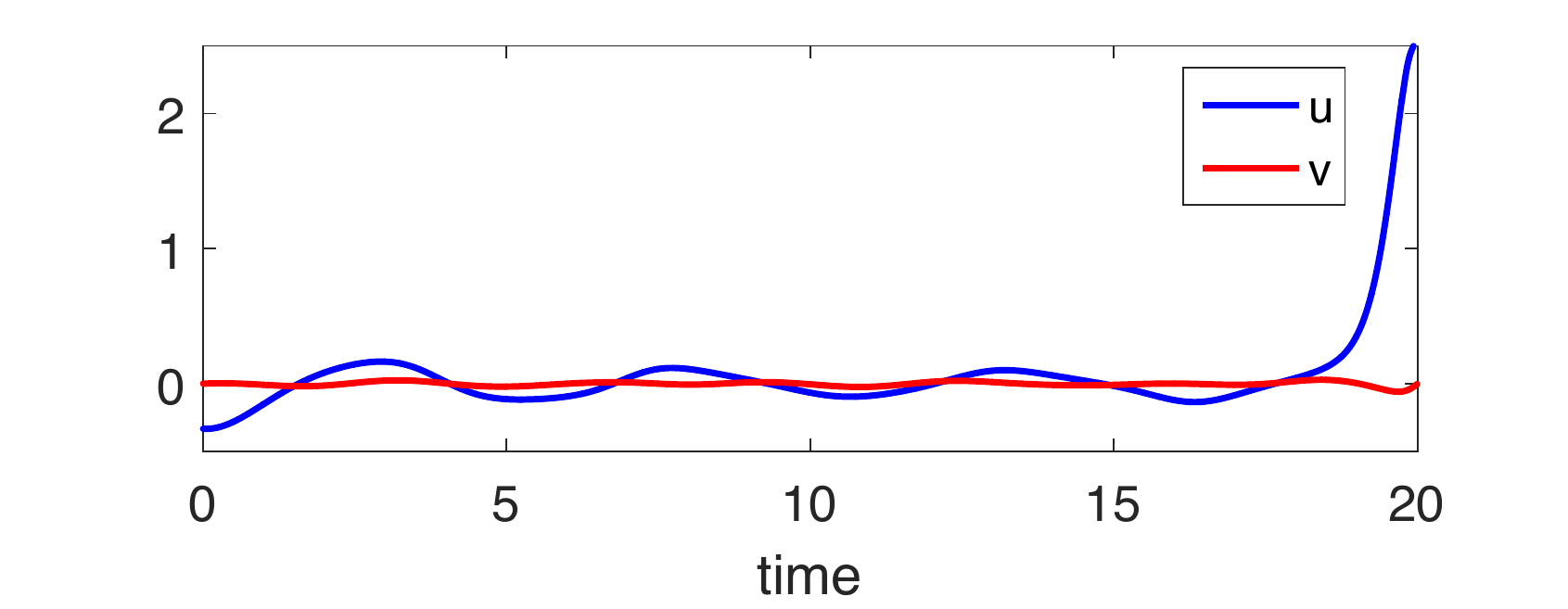}
	\label{fig:ensemblecontrol}}
	\subfigure[The final x-components of the ensemble]{\includegraphics[width=\columnwidth]{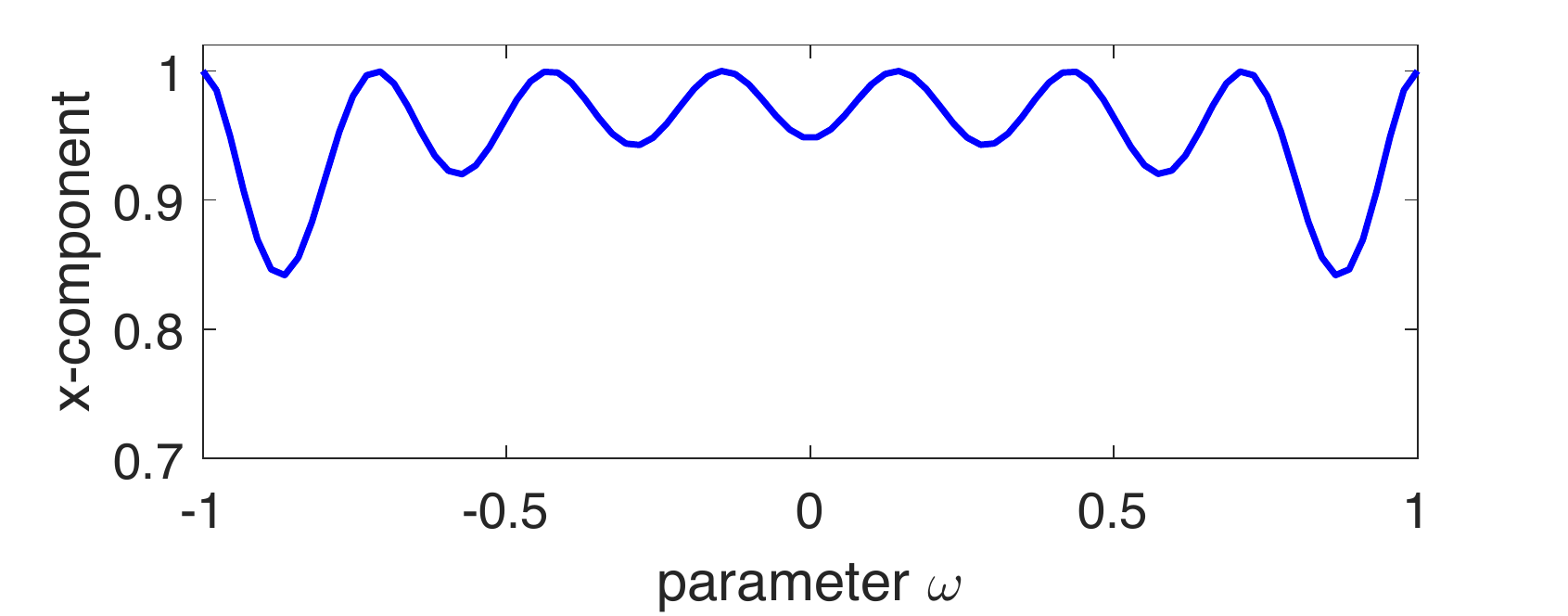} 
	\label{fig:ensemblefinalx}}
	\caption{\subref{fig:ensemblecontrol} The optimal ensemble control that steers an ensemble of Bloch systems with $\w\in[-1,1]$ from $X_0(\w) = (0,0,1)^T$ to a neighborhood of $X_d(\w)=(1,0,0)^T$. The weighted matrix $R=I_3$. \subref{fig:ensemblefinalx} The final states $X(20,\w)$ for 81 spin systems (blue) with their frequencies uniformly spaced within $[-1,1]$ following the control displayed in \subref{fig:ensemblecontrol}.}
	\label{fig:blochensemble}
\end{figure}

We apply the iterative method described in Section \ref{sec:iterative} to design a broadband excitation ($\pi/2$) pulse that minimizes the cost functional as in \eqref{eq:J}, where $u=(u_1,u_2)^T$ and the terminal cost $J_E=\int_{-1}^1 [X(t_f,\w)-X_d(\w)]^T [X(t_f,\w)-X_d(\w)]\, d\w$ with $t_f=20$. Figure \ref{fig:ensemblecontrol} shows the derived broadband pulse with $R=I_{3}$. The performance, i.e., the $x$-components of the final states, are displayed in Figure \ref{fig:ensemblefinalx}. The iterative algorithm converges in 207 iterations given the stopping criterion $\|X^{(k+1)} - X^{(k)}\| < 10^{-4}$.
\end{example}

\begin{figure}[t]
  	\centering
  	\subfigure[The optimal control law]{\includegraphics[width=1\columnwidth]{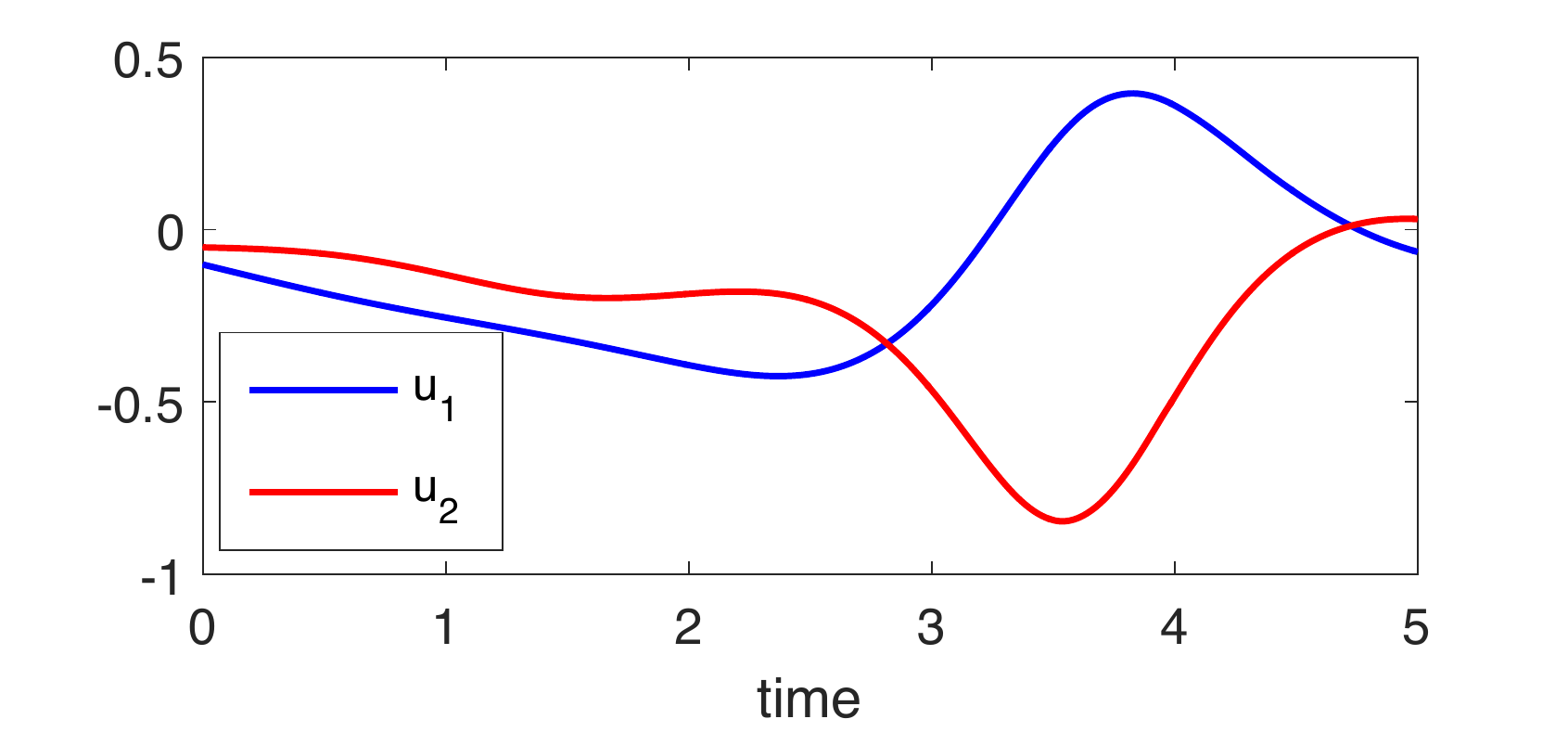} \label{fig:control}}
  	\subfigure[Coherence transfer]{\includegraphics[width=1\columnwidth]{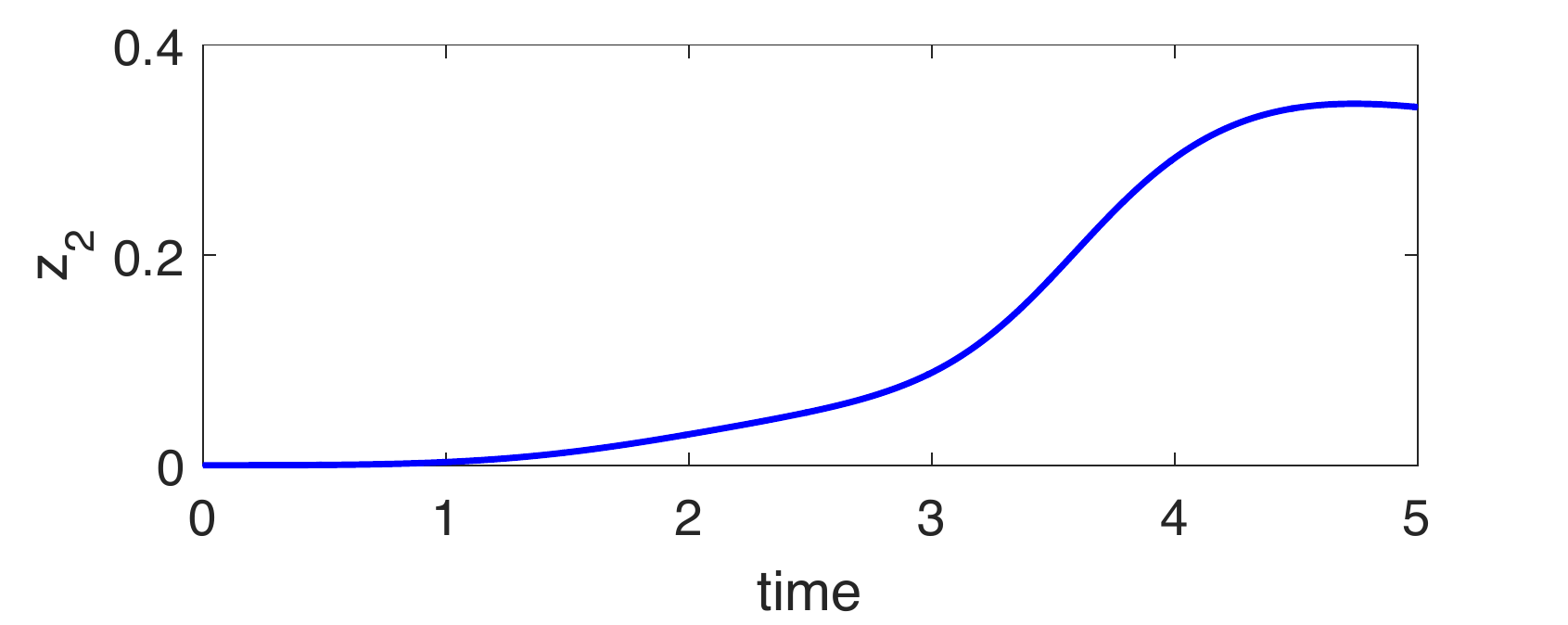}  \label{fig:path}}
	\subfigure[The convergence behavior]{\includegraphics[width=1\columnwidth]{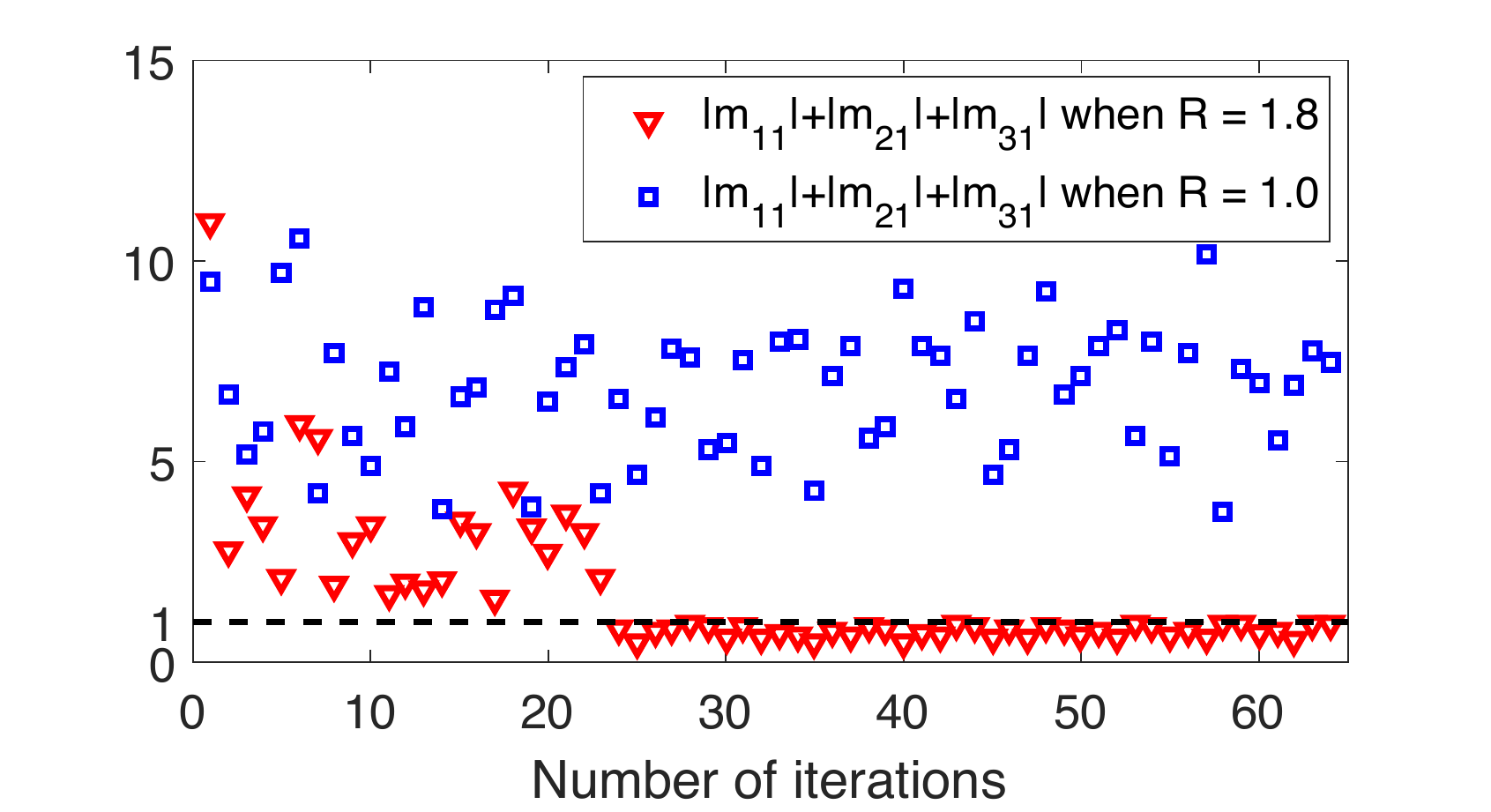} \label{fig:CompareM}}
  	\caption{\subref{fig:control} The optimal control that maximizes coherence transfer for the case in which the weighted matrix $R=1.8I_6$, $\omega_1 = \omega_2 = 0.5$, $\xi_c = 0.8$, and $\xi_a=1$. \subref{fig:path} The resulting $z_2$ trajectory following the optimal control in \subref{fig:control}. \subref{fig:CompareM} The convergence behavior corresponding to the cases of $R=1.8I_6$ and $R=I_6$, where in the former case the convergence criterion in \eqref{eq:criterion} is satisfied after 24 iterations. }
  \label{fig:twobloch}
  \end{figure}

\begin{example}[Coherence Transfer of Two-Spin Systems]
\label{ex:coherence}
\rm Transferring coherence between two spins is a fundamental step to multidimensional NMR spectroscopy. In protein NMR spectroscopy, large transverse relaxation rates can cause degraded sensitivity and thereby limit the size of macromolecules available for study \cite{BBCROP}. Coherence transfer between a spin pair can be modeled as an optimal control problem involving a bilinear system that described the spin dynamics of the form $\dot{X}(t)=A(X(t),u)X(t)$, where $X = (z_1,y_1,x_1,x_2,y_2,z_2)^T$ denotes the state (expectation values of the spin operators \cite{Li_JCP09}) of the two-spin system with the natural frequencies $\w_1$ and $\w_2$, respectively; $u = (u_1,u_2)^T$ is the control pulse, and 
\begin{align}
	\label{eq:twospin}
 	A = \begin{bmatrix}
 	0 & -u_1 & u_2 & 0 & 0 & 0\\
 	u_1 & -\xi_a & \omega_1 & J & -\xi_c & 0\\
 	-u_2 & -\omega_1 & -\xi_a & -\xi_c & J & 0 \\
 	0 & J & -\xi_c & -\xi_a& -\omega_2 & -u_2 \\
 	0 & -\xi_c & J & \omega_2 & -\xi_a & u_1\\
 	0 & 0 & 0 & u_2 & -u_1 & 0
 	\end{bmatrix}.
\end{align}
The parameters $J$, $\xi_a$ and $\xi_c$ 
represent the scalar coupling constant, the transverse autocorrelated relaxation rate, and the cross-correlation relaxation rate.

%
%

The objective is to design a relaxation-optimized pulse $u$ that maximizes the coherence transfer, i.e., maximizes the final value of $x_2$ at time $t_f$, free to vary as a decision variable, from the initial state $X_0 = [1,0,0,0,0,0]^T$. In particular, we consider the cost functional $\mathcal{J}=\int_0^{t_f} \frac{1}{2}u^T(t)Ru(t)dt+\|X(tf)-X_d\|^2$, where $X_d = (0,0,0,0,0,1)^T$, with the tradeoff between coherence transfer and pulse energy. Applying the iterative algorithm with $R=1.8I_6$ results in the convergent optimal control and optimal trajectory of $x_2(t)$ shown in Figure \ref{fig:control} and \ref{fig:path}, respectively, where the final time is $t_f=5$ and the parameters are $J=0.5$, $\xi_a=1$, $\xi_c=0.8$, and $\w_1=\w_2=0.5$. The algorithm converges in 64 iterations given the stopping criterion $\|X(t_f)-X_d\|\leq 10^{-3}$. The maximum transfer is $x_2=0.3425$ within the time horizon $t\in [0,5]$.

Moreover, we use this example to illustrate the convergence condition \eqref{eq:criterion} by comparing the convergence behavior between difference choices of the control weight matrix. Here, we use $R=1.8I_6$ and $R=I_6$. In the latter case, the iterative algorithm fails to be convergent since the condition \eqref{eq:criterion} is violated as shown in Figure \ref{fig:CompareM}, while in the former case, the condition \eqref{eq:criterion} is satisfied for $k\geq 24$. 
This implies that $(X^{(k)},K^{(k)},s^{(k)})$ enters the invariant set containing the fixed point after 24 iterations.

\end{example}

\section{Conclusion}
In this paper, we study optimal control problems with free-endpoint conditions involving an inhomogeneous bilinear ensemble system with a parameter-dependent drift. We develop an iterative method to solve this class of analytically intractable problems. The central idea is to represent the bilinear system at each iteration as a time-varying linear system and then solve it in an iterative manner. We show the convergence of the iterative procedure by the use of the fixed-point theorem, and illustrate that the magnitude of the eigenvalues of the control weight matrix $R$ plays the key role in the determination of the convergence. Our method is directly applicable to study stochastic optimal ensemble control problems driven by additive noise processes such as the Gaussian and Poisson counting processes.

\section{Appendix} 
\subsection{The $\b$ Coefficients in the Proof of Theorem 1}
\label{appd:betas}
Let $\Phi^{(k)}$ be the transition matrix associated with the system $\dot{x} = -\big[\tilde{A}^{(k)} - \tilde{O}^{(k)}K^{(k+1)} \big]^T x$. The time-varying coefficients $\beta$'s are as follows:
\begin{align*}
\beta_1 = &  \|\Phi^{(k)}(\s,t)\| \|s^{(k)}(\s)\|, \\
\beta_2 = & \|\Phi^{(k)}(\s,t)\| \|K^{(k+1)}(\s)\| \|s^{(k)}(\s)\|, \\
\beta_3 = & \|\Phi^{(k)}(\s,t)\| \| \tilde{O}^{(k-1)} (\s)\| \|s^{(k)}(\s)\| + \|\Phi^{(k)}(\s,t)\| \| g\|, \\
\beta_4 = &  \|\Phi^{(k)}(\s,t)\| ( \|K^{(k)}(\s)\|+\|K^{(k+1)}(\s)\| ) \|\Phi^{(k)}(\s,t)\| , \\
\beta_5 =  &  \|\Phi^{(k)}(\s,t)\|\|K^{(k)}(\s)\|\|K^{(k+1)}(\s)\| \|\Phi^{(k)}(\s,t)\| , \\
\beta_6 = &  \|\Phi^{(k)}(t,\s)\| \|x^{(k)}(\s)\|, \\
\beta_7 = & \|\Phi^{(k)}(t,\s)\| \|\tilde{O}^{(k-1)} (\s)\| \|x^{(k)}(\s)\|, \\
\beta_8 = &  \|\Phi^{(k)}(t,\s)\| \Big( \|K^{(k+1)}(\s)\|\|x^{(k)}(\s)\| + \|s^{(k+1)}(\s)\| \Big) , \\
\beta_9 = &  \|\Phi^{(k)}(t,\s)\| \|\tilde{O}^{(k-1)} (\s)\|\|x^{(k)}(\s)\|.
\end{align*}

\subsection{The entries of the matrix $M$ for the contraction mapping}
\label{appd:M}
Let $\delta=\sqrt{\sum_{i=1}^n \|P_i \| ^2}$ and $\zeta=\sqrt{\sum_{i,j=1}^n \|Q_{ij} \| ^2}$,which appeared in \eqref{eq:dA} and \eqref{eq:dO}, then the entries of the matrix $M$ defined in \eqref{eq:contraction} satisfy the following properties
\begin{align*}
m_{11} & \varpropto (\beta_6 + \beta_4 \beta_7 + \b_1 \b_9 + \b_3 \b_4 \b_9) \Big( (\delta+\|x^{(k-1)} \|\zeta)\|K^{(k)} \|   \\
& + \zeta (\|K^{(k)}\|\|x^{(k)}\|+\|s^{(k)}\|) \Big) + \zeta \Big( \|x^{(k)}\|+\|x^{(k-1)}\| \Big) \cdot  \\
& \cdot (\beta_8 + \beta_5 \beta_7 + \b_2 \b_9 + \b_3 \b_5 \b_9), \\
m_{12} & \varpropto (\beta_6 + \beta_4 \beta_7 + \b_1 \b_9 + \b_3 \b_4 \b_9)(\delta+\|x^{(k-1)} \|\zeta )\|x^{(k-1)} \|,\\
m_{13} & \varpropto (\beta_6 + \beta_4 \beta_7 + \b_1 \b_9 + \b_3 \b_4 \b_9)(\delta+\|x^{(k-1)} \|\zeta ),\\ 
m_{21} & \varpropto \beta_4 \Big( (\delta+\|x^{(k-1)} \|\zeta )\|K^{(k)} \| + \zeta \|K^{(k)}\|\|x^{(k)}\|  \\
& + \zeta \|S^{(k)}\|\Big) + \beta_5 \zeta \Big( \|x^{(k)}\|+\|x^{(k-1)}\| \Big), \\
m_{22} & \varpropto \beta_4 (\delta+\|x^{(k-1)} \|\zeta )\|x^{(k-1)} \|, \\
m_{23} & \varpropto \beta_4 (\delta+\|x^{(k-1)} \|\zeta ), \\
m_{31} & \varpropto (\beta_1 + \beta_3 \beta_4) \Big( (\delta+\|x^{(k-1)} \|\zeta )\|K^{(k)} \| + \zeta \|K^{(k)}\|\|x^{(k)}\|  \\
& + \zeta \|s^{(k)})\| \Big) + \zeta (\beta_2 + \beta_3 \beta_5)  \Big( \|x^{(k)}\|+\|x^{(k-1)}\| \Big), \\
m_{32} & \varpropto (\beta_1 + \beta_3 \beta_4) (\delta+\|x^{(k-1)} \|\zeta )\|x^{(k-1)} \|, \\
m_{33} & \varpropto (\beta_1 + \beta_3 \beta_4) (\delta+\|x^{(k-1)} \|\zeta ).
\end{align*}

\subsection{The proof for Theorem \ref{thm:globalsufficient}}
\label{appd:proof}
{\it Proof:}
Because the partial derivatives of $V^{(k)}$ in \eqref{eq:Vfunc} with respect to $t$ and $x$ can be expressed in terms of $x^{(k)}$, $K^{(k)}$, and $s^{(k)}$, which are convergent according to Theorem \ref{thm:convergence}, the sequences $\Big\{\dfrac{\partial V^{(k)}}{\partial t}(t,x^{(k)}(t))\Big\}_k$ and $\Big\{ \dfrac{\partial V^{(k)}}{\partial x}(t,x^{(k)}(t))\Big\}_k$ are convergent, namely, $\dfrac{\partial V^{(k)}}{\partial t}(t,x^{(k)}(t)) \rightarrow V_t(t,x^*)$ and $\dfrac{\partial V^{(k)}}{\partial x}(t,x^{(k)}(t))\rightarrow V_x(t,x^*)$, where
\begin{align}
	\nonumber
	V_t(t,x^*) &= \frac{1}{2}(x^*)^T \Big[- (\tilde{A}^*)^TK^* - K^*\tilde{A}^* + K^*\tilde{O}^*K^*\Big] x^* \\
	\nonumber
	&+ (x^*)^T \Big[ -(\tilde{A}^*)^T - K^*\tilde{O}^* \Big]s^* + \dfrac{1}{2} (s^*)^T \tilde{O}^* s^* ,\\
	\label{eq:p*}
	V_x(t,x^*) &= K^*x^* +s^* = p^*,
\end{align}
in which $\tilde{A}^{(k)}\rightarrow\tilde{A}^*$, $\tilde{O}^{(k)}=\tilde{B}^{(k)} R^{-1}(\tilde{B}^{(k)})^T\rightarrow\tilde{B}^*R^{-1}(\tilde{B}^*)^T=\tilde{O}^*$, and $p^{(k)}\rightarrow p^*$ follow from \eqref{eq:Ak}, \eqref{eq:Bk}, and \eqref{eq:p^k}, respectively.
Because $\dfrac{\partial V^{(k)}}{\partial t}(t,x^{(k)}(t))$ and $\dfrac{\partial V^{(k)}}{\partial x}(t,x^{(k)}(t))$ are dominated by $\xi(t)$ and $\eta(x)=(\eta_1,\ldots,\eta_n)^T$, respectively, by the Lebesgue Dominated Convergence theorem, 
we have
\begin{align}
	V^* (t,x) & = \lim_{k\rightarrow \infty} V^{(k)}(t,x) = \lim_{k\rightarrow \infty}\int_0^t \dfrac{\partial V^{(k)}}{\partial \s}(\s,x) d\s\nonumber \\ 
	\label{eq:V*1}
	& =\int_0^t \lim_{k\rightarrow \infty}\dfrac{\partial V^{(k)}}{\partial \s}(\s,x) d\s = \int_0^t V_t (\s,x) d\s, \\
	V^* (t,x) & = \lim_{k\rightarrow \infty} V^{(k)}(t,x) = \lim_{k\rightarrow \infty}\int_{x(0)}^{x(t)} \dfrac{\partial V^{(k)}}{\partial x}(t,x) dx \nonumber \\
	\label{eq:V*2}	
	& =\int_{x(0)}^{x(t)} \lim_{k\rightarrow \infty}\dfrac{\partial V^{(k)}}{\partial x}(t,x) dx = \int_{x(0)}^{x(t)} V_x (t,x) dx.
\end{align}
Because $V^*$ is continuously differentiable with respect to both $t$ and $x$, we obtain from \eqref{eq:V*1} and \eqref{eq:V*2} the partial derivatives
\begin{align} 
	\label{eq:direvative}
	\dfrac{\partial V^*}{\partial t}(t,x) = V_t (t,x),\quad \dfrac{\partial V^*}{\partial x}(t,x) = V_x (t,x).
\end{align}
In addition, due to the convergence of the iterative procedure, \eqref{eq:fopde} is convergent to
\begin{align*}
	V_t(t,x^*) &+ V_x(t,x^*)^T(\tilde{A}^*x^* + g) - \frac{1}{2} V_x(t,x^*)^T \tilde{O}^* V_x(t,x^*)\equiv 0.
\end{align*}
which, by employing \eqref{eq:direvative} and \eqref{eq:p*}, can be rewritten as,
\begin{align*}
	V_t(t,x^*) +  (p^*)^T (\tilde{A}^*x^*-\tilde{O}^*p^*+g)+\dfrac{1}{2}(p^*)^T \tilde{O}^*p^*\equiv 0,
\end{align*}
with the boundary condition $V^*(t_f,x^*)=0$. Because the convergent solution pair $(x^*,p^*)$ satisfies the necessary condition, the above equation is equivalent to, by \eqref{eq:optimalu}, \eqref{eq:state1} and \eqref{eq:costate1},
\begin{align} 
	&\dfrac{\partial V^*}{\partial t}(t,x^*) +  \dfrac{\partial V^*}{\partial x}(t,x^*)^T \Big[Ax^*+\Big(B+\sum_{i=1}^n x^*_i N_i\Big) u^* +g\Big]\nonumber \\
	\label{eq:HJBu}
	&+\dfrac{1}{2}(u^*)^TRu^* \equiv 0.
\end{align}
Since $V^*$ is differentiable, according to the principle of dynamic programming, the quantity on the left-hand side in \eqref{eq:HJBu} is non-negative for every control $u$ in the admissible control set $\mathcal{U}\subseteq PC([0,t_f];\, \mathbb{R}^m)$. It follows that $(V^*,u^*)$ is a solution to the HJB equation of the original Problem \eqref{eq:oc1}, that is,
\begin{align} 
	&\dfrac{\partial V}{\partial t}(t,x) + \min_{u\in\mathcal{U}} \Big\{ \dfrac{\partial V}{\partial x}(t,x)^T \Big[Ax+\Big(B+\sum_{i=1}^n x_i N_i\Big) u+g\Big]\nonumber \\
	\label{eq:HJBbilinear}
	& +\dfrac{1}{2} u^TRu  \Big\} \equiv 0,
\end{align}
with the boundary condition $V(t_f,x)=0$. Furthermore, the optimal control $u^*$ is global since the minimization in \eqref{eq:HJBbilinear} is over a convex (quadratic) function in $u$. This proves that the convergent solution $u^*$ of the iterative method 
is a global minimizer to the original Problem \eqref{eq:oc1}. \hfill$\Box$

\bibliographystyle{plain}        
\bibliography{bilinear_ensemble}           



\appendix
\end{document}